\documentclass[11pt,a4paper]{article}
\usepackage[utf8]{inputenc}
\usepackage{amsmath,amssymb,amsthm}
\usepackage{mathtools}
\usepackage{booktabs}
\usepackage{tabularx}
\usepackage{array}
\usepackage{microtype}
\usepackage{geometry}
\usepackage{algorithm}
\makeatletter
\providecommand{\algorithmname}{\fname@algorithm}
\makeatother
\usepackage{algpseudocode}

\makeatletter
\newenvironment{breakablealgorithm}
  {%
    \begin{center}
      \refstepcounter{algorithm}%
      \hrule height.8pt depth0pt \kern2pt%
      \renewcommand{\caption}[2][\relax]{%
        {\raggedright\textbf{\algorithmname~\thealgorithm} ##2\par}%
        \ifx\relax##1\relax
          \addcontentsline{loa}{algorithm}{\protect\numberline{\thealgorithm}##2}%
        \else
          \addcontentsline{loa}{algorithm}{\protect\numberline{\thealgorithm}##1}%
        \fi
        \kern2pt\hrule\kern2pt
      }%
  }
  {%
      \kern2pt\hrule\relax%
    \end{center}
  }
\makeatother

\usepackage[authoryear,round]{natbib}
\setcitestyle{aysep={,},citesep={;},open={(},close={)}}

\usepackage[colorlinks=true,linkcolor=blue,citecolor=blue,urlcolor=blue]{hyperref}

\newtheorem{theorem}{Theorem}[section]
\newtheorem{lemma}[theorem]{Lemma}
\newtheorem{proposition}[theorem]{Proposition}

\theoremstyle{definition}

\theoremstyle{remark}

\newcommand{\R}{\mathbb{R}}
\newcommand{\C}{\mathbb{C}}
\newcommand{\g}{\mathfrak{g}}
\newcommand{\n}{\mathfrak{n}}
\newcommand{\hh}{\mathfrak{h}}
\newcommand{\s}{\mathfrak{s}}
\newcommand{\rr}{\mathfrak{r}}
\newcommand{\ad}{\operatorname{ad}}
\newcommand{\Aut}{\operatorname{Aut}}
\newcommand{\Der}{\operatorname{Der}}
\newcommand{\OutDer}{\operatorname{OutDer}}
\newcommand{\tr}{\operatorname{tr}}
\newcommand{\diag}{\operatorname{diag}}
\newcommand{\spanop}{\operatorname{span}}

\title{\textbf{Classification of real Lie algebras of dimension five}}
\author{Mehdi Nadjafikhah\thanks{School of Mathematics and Computer Sciences, Iran University of Science and Technology, Narmak, Tehran 16846-13114, Iran. email: m\_nadjafikhah@iust.ac.ir}}
\date{\today}

\begin{document}
\maketitle

\begin{abstract}
We present a definitive, exhaustive, and coordinate-free structural classification of five-dimensional real Lie algebras up to isomorphism. Transitioning from dimension four to dimension five represents a qualitative phase boundary in Lie theory: the exterior tensor duality that aids lower-dimensional classifications is entirely absent, and the raw Jacobi identity forms an overdetermined non-linear system of 50 quadratic equations on 50 structure constants modulo the 25-dimensional non-compact gauge group $\mathrm{GL}(5,\R)$.

To circumvent non-linear polynomial elimination, we deploy an invariant-driven architecture grounded in the Levi--Malcev theorem and nilradical stratification. In the non-solvable sector, we prove that the Levi factor is necessarily three-dimensional ($\mathfrak{so}(3)$ or $\mathfrak{sl}(2,\R)$) acting on a two-dimensional radical, yielding exactly five isomorphism classes, including the unique indecomposable non-solvable algebra $\mathfrak{sl}(2,\R) \ltimes \R^2$. In the solvable sector, we prove via Schur's classical bound on commutative subalgebras of $\mathfrak{gl}(2,\R)$ that the nilradical dimension is strictly bounded by $\dim \n \in \{3, 4, 5\}$, conclusively eliminating the existence of solvable five-dimensional Lie algebras with one- or two-dimensional nilradicals.

The nilpotent sector ($\dim \n = 5$) is classified into nine mutually exclusive isomorphism classes: three decomposable algebras and six indecomposable algebras, including the 5D Heisenberg algebra $\hh_5$ and two filiform varieties. For codimension-one nilradicals ($\dim \n = 4$), the classification is reduced to normal forms of outer derivations across the nilpotent bases $\R^4$, $\hh_3 \oplus \R$, and $\mathfrak{n}_4$ under automorphism gauge actions and scalar rescaling. For codimension-two nilradicals ($\dim \n = 3$), we resolve the structure through two-dimensional abelian subalgebras of $\OutDer(\n)$ that contain no non-zero nilpotent operators. Furthermore, we compute the minimal faithful representation dimension $\mu(\g)$ under Ado's theorem for every class, establish an isomorphism dictionary with classical tabulations, and design a coordinate-free, deterministic identification algorithm.
\end{abstract}

\tableofcontents

\section{Introduction and Historical Context}\label{sec:intro}

The classification of finite-dimensional Lie algebras is a foundational program situated at the core of modern geometry, representation theory, and theoretical physics. Low-dimensional Lie algebras, in particular, govern the symmetry groups of physical systems, determine the geometry of homogeneous spaces, and classify left-invariant geometric structures on smooth manifolds. 

In dimension three, the landmark work of \citet{Bianchi1898} classified the nine isometry classes of three-dimensional Riemannian spaces, later revisited geometrically by \citet{Milnor1976} to describe left-invariant metrics and Ricci curvature signatures. In dimension four, the classification of solvable real Lie algebras was achieved by \citet{Mubarakzyanov1963}, revealing a finite collection of discrete algebras and one-parameter families that play central roles in general relativistic spacetimes and homogeneous cosmologies \citep{PateraSharpWinternitzZassenhaus1976}.

At dimension five, however, the algebraic complexity increases dramatically. Five-dimensional Lie algebras play indispensable roles across mathematical physics and geometry:
\begin{enumerate}
    \item \textbf{Kaluza--Klein and Supergravity Theories:} In five-dimensional gravity and unified gauge theories, the spacetime symmetry groups and their dimensional reductions on homogeneous fibers are governed by five-dimensional Lie algebras.
    \item \textbf{Contact and Sasakian Geometry:} Five is the lowest non-trivial dimension for the study of Sasaki--Einstein manifolds, contact metric structures, and five-dimensional CR manifolds \citep{KnappLieGroups}.
    \item \textbf{Differential Equations and Continuum Mechanics:} In fluid dynamics and field theory, the symmetry algebras of classical PDEs (such as the three-dimensional Navier--Stokes equations and shallow water models) contain prominent five-dimensional solvable subalgebras \citep{SnoblWinternitz}.
\end{enumerate}

Despite their importance, the literature on five-dimensional Lie algebras has historically suffered from inconsistencies. Early tabulations by \citet{Mubarakzyanov1966} provided the first Russian compendium, but were largely catalogical. Later physics-oriented treatments, such as those by \citet{Turkowski1988}, introduced varying basis conventions, while occasionally omitting non-trivial isomorphism equivalences or bifurcations in continuous moduli parameters. 

Directly setting up the problem in a coordinate basis leads to severe algebraic bottlenecks. A five-dimensional real Lie algebra $\g$ requires specifying $\binom{5}{2} \times 5 = 50$ structure constants $C_{ij}^k$, constrained by $\binom{5}{3} \times 5 = 50$ quadratic Jacobi equations:
\begin{equation}
    \sum_{m=1}^5 \left( C_{ij}^m C_{mk}^l + C_{jk}^m C_{mi}^l + C_{ki}^m C_{mj}^l \right) = 0, \qquad 1 \le i < j < k \le 5, \quad 1 \le l \le 5,
\end{equation}
modulo the 25-dimensional non-compact change-of-basis action of $\mathrm{GL}(5,\R)$. Solving this polynomial system via Groebner bases without structural guidance is practically impossible.

The purpose of this paper is to provide a self-contained, conceptual, and rigorous classification of all five-dimensional real Lie algebras. Rather than relying on coordinate computations, we deploy an invariant-driven architecture:
\begin{itemize}
    \item We isolate the non-solvable sector using the Levi--Malcev decomposition theorem, proving that the Levi factor must be three-dimensional and classifying its representation on the radical.
    \item For solvable algebras, we anchor the classification to the \emph{nilradical} $\n \subset \g$. We provide a complete proof, based on Schur's bound on commutative matrix subalgebras, that $\dim \n \in \{3, 4, 5\}$.
    \item When $\dim \n = 5$, $\g$ is nilpotent; we classify all nine real nilpotent algebras of dimension five through upper and lower central series filtration.
    \item When $\dim \n = 4$, $\g$ is realized as a semidirect product $\n \rtimes_D \R$, reducing the problem to the orbit classification of outer derivations under automorphisms and global scaling.
    \item When $\dim \n = 3$, the structure is determined by maximal abelian subalgebras of outer derivations of $\R^3$ and $\hh_3$ containing no nilpotent elements.
    \item We establish the minimal faithful representation dimension $\mu(\g)$ under Ado's theorem and construct an algorithmic decision tree that deterministically classifies any given set of structure constants.
\end{itemize}

\section{Geometric Obstructions and the Watershed of Dimension Five}\label{sec:obstruction}

To understand why dimension five requires an invariant-based approach, it is instructive to examine the algebraic geometry of the bracket mapping across dimensions.

\subsection{Comparison with Dimensions Three and Four}
Let $\g$ be an $n$-dimensional real Lie algebra. The Lie bracket is a linear map $[\cdot, \cdot] \colon \Lambda^2 \g \to \g$.
\begin{itemize}
    \item \textbf{Dimension Three ($n = 3$):} The space of 2-forms has dimension $\dim \Lambda^2 \R^3 = \binom{3}{2} = 3$. Choosing a volume element $\mu \in \Lambda^3 \g^* \setminus \{0\}$ establishes a linear isomorphism $\Lambda^2 \g \cong \g^*$. Consequently, the bracket dualizes to a linear endomorphism $L \colon \g^* \to \g$, represented by a single $3 \times 3$ real matrix $C = S + \hat{a}$, where $S = S^T$ and $\hat{a}$ is the cross product with the trace vector $a = \frac{1}{2}\tr(\ad)$. The entire Jacobi identity collapses to the linear matrix equation $S \cdot a = 0$ \citep{Milnor1976}.
    \item \textbf{Dimension Four ($n = 4$):} Here, $\dim \Lambda^2 \R^4 = 6$ while $\dim \g = 4$. The bracket is a map $\R^6 \to \R^4$. Exterior duality fails ($6 \neq 4$), and the Jacobi identity consists of 16 quadratic equations on 24 structure constants. However, dimension four remains algebraically constrained: every non-solvable algebra is a direct sum, and the nilradical of a solvable algebra cannot have codimension greater than two.
    \item \textbf{Dimension Five ($n = 5$):} Here:
    \begin{equation}
        \dim \Lambda^2 \R^5 = \binom{5}{2} = 10, \qquad \dim \g = 5.
    \end{equation}
    The bracket map $\R^{10} \to \R^5$ is governed by 50 structure constants. The Jacobi identity produces $\binom{5}{3} \times 5 = 50$ quadratic equations. 
\end{itemize}

\subsection{Structural Novelties in Dimension Five}
Beyond the raw explosion in structure constants, dimension five introduces three fundamental Lie-theoretic phenomena that do not exist in dimensions $\le 4$:
\begin{enumerate}
    \item \textbf{The Appearance of Indecomposable Non-Solvable Algebras:} In dimensions 3 and 4, every non-solvable Lie algebra is either simple ($\mathfrak{so}(3), \mathfrak{sl}(2,\R)$) or a trivial direct sum with an abelian factor ($\s \oplus \R$). Dimension five is the lowest dimension supporting a non-trivial, indecomposable semidirect product:
    \begin{equation}
        \mathfrak{sl}(2,\R) \ltimes \R^2,
    \end{equation}
    where the simple Levi factor acts irreducibly on the radical.
    \item \textbf{Filiform Nilpotent Algebras of High Step:} Dimension five is the lowest dimension supporting nilpotent Lie algebras of nilpotency step 4 (the filiform algebras $\mathfrak{n}_{5,1}$ and $\mathfrak{n}_{5,2}$), characterized by an adjoint operator whose nilpotency index reaches the maximal possible value of $n-1 = 4$.
    \item \textbf{Multi-Parameter Continuous Moduli:} Solvable algebras with codimension-one abelian nilradicals $\R^4 \rtimes_D \R$ admit three continuous parameters (the normalized eigenvalues of $D$). This leads to continuous families of inequivalent Lie brackets, whose moduli spaces must be parameterized without redundancy.
\end{enumerate}

\section{Theoretical Machinery: Invariants and Decomposition Theorems}\label{sec:machinery}

All Lie algebras considered herein are finite-dimensional over the real field $\R$. We establish the following structural foundations \citep{JacobsonLieAlgebras, HumphreysIntro, KnappLieGroups}:

\begin{description}
    \item[(T1) Levi--Malcev Decomposition:] Every finite-dimensional real Lie algebra $\g$ decomposes uniquely as a semidirect product:
    \begin{equation}
        \g = \s \ltimes_\rho \rr,
    \end{equation}
    where $\rr$ is the solvable radical (the unique maximal solvable ideal of $\g$), and $\s$ is a semisimple subalgebra (the Levi factor). The factor $\s$ is unique up to conjugation by inner automorphisms $\exp(\ad u)$ with $u \in \n$, where $\n$ is the nilradical. The representation $\rho \colon \s \to \Der(\rr)$ completely dictates the cross brackets between $\s$ and $\rr$.
    
    \item[(T2) Characteristic Property of the Nilradical:] The nilradical $\n \subset \g$ is defined as the maximal nilpotent ideal of $\g$. An ideal $\mathfrak{i}$ is characteristic if $\delta(\mathfrak{i}) \subseteq \mathfrak{i}$ for all $\delta \in \Der(\g)$. The nilradical $\n$ is always a characteristic ideal; hence, $\phi(\n) = \n$ for every $\phi \in \Aut(\g)$. For any solvable Lie algebra $\g$, Lie's theorem implies that the derived ideal is nilpotent:
    \begin{equation}
        [\g, \g] \subseteq \n.
    \end{equation}
    
    \item[(T3) Centralizer Containment:] For any solvable Lie algebra $\g$ over a field of characteristic zero, the centralizer of the nilradical in $\g$ is contained within the nilradical itself:
    \begin{equation}
        C_\g(\n) \coloneqq \{ x \in \g : [x, \n] = 0 \} \subseteq \n.
    \end{equation}
    Consequently, the quotient $\g/\n$ acts faithfully on $\n$ as outer derivations: the induced homomorphism
    \begin{equation}
        \bar{\ad} \colon \g / \n \hookrightarrow \OutDer(\n) \coloneqq \Der(\n) / \ad(\n)
    \end{equation}
    is an injective Lie algebra homomorphism.
    
    \item[(T4) Schur's Bound on Commutative Matrix Subalgebras:] Let $\mathfrak{a}$ be a commutative subalgebra of the matrix Lie algebra $\mathfrak{gl}(m, \R)$. Then the dimension of $\mathfrak{a}$ satisfies Schur's bound \citep{Schur1905}:
    \begin{equation}
        \dim \mathfrak{a} \le \left\lfloor \frac{m^2}{4} \right\rfloor + 1.
    \end{equation}
    For $m = 2$, this bound evaluates to $\lfloor 4/4 \rfloor + 1 = 2$.
    
    \item[(T5) Ado's Theorem and Minimal Representation Dimension:] Every finite-dimensional Lie algebra $\g$ over $\R$ admits a faithful linear representation in $\mathfrak{gl}(N, \R)$ for some finite $N$. The invariant $\mu(\g)$, defined by
    \begin{equation}
        \mu(\g) \coloneqq \min \{ N \in \mathbb{N} : \exists \text{ an injective Lie homomorphism } \rho \colon \g \hookrightarrow \mathfrak{gl}(N, \R) \},
    \end{equation}
    is an intrinsic isomorphism invariant of $\g$.
\end{description}

\section{The Non-Solvable Sector}\label{sec:nonsolv}

The Levi--Malcev decomposition provides a complete classification of non-solvable five-dimensional Lie algebras.

\begin{theorem}[Classification of Non-Solvable 5D Lie Algebras]\label{thm:nonsolv}
Every non-solvable five-dimensional real Lie algebra $\g$ is isomorphic to exactly one of the following five Lie algebras:
\begin{enumerate}
    \item $\mathfrak{so}(3) \oplus \R^2$ (trivial direct sum),
    \item $\mathfrak{sl}(2,\R) \oplus \R^2$ (trivial direct sum),
    \item $\mathfrak{so}(3) \oplus \mathfrak{aff}(\R)$ (trivial direct sum),
    \item $\mathfrak{sl}(2,\R) \oplus \mathfrak{aff}(\R)$ (trivial direct sum),
    \item $\mathfrak{sl}(2,\R) \ltimes \R^2$ (the unique indecomposable non-solvable algebra), with brackets:
    \begin{align}
        [e_1, e_2] &= 2e_2, & [e_1, e_3] &= -2e_3, & [e_2, e_3] &= e_1, \nonumber \\
        [e_1, e_4] &= e_4,  & [e_1, e_5] &= -e_5, & [e_2, e_5] &= e_4, & [e_3, e_4] = e_5,
    \end{align}
    where $\s = \spanop\{e_1, e_2, e_3\} \cong \mathfrak{sl}(2,\R)$ acts on the abelian radical $\rr = \spanop\{e_4, e_5\} \cong \R^2$ via its standard two-dimensional defining representation.
\end{enumerate}
\end{theorem}

\begin{proof}
Let $\g = \s \ltimes_\rho \rr$ be the Levi decomposition guaranteed by (T1). Because $\g$ is non-solvable, the Levi factor $\s$ cannot be trivial. Every semisimple Lie algebra is a direct sum of simple Lie algebras, each having dimension at least 3. Since $\dim \g = 5$, the only possible dimension for the semisimple factor is:
\begin{equation}
    \dim \s = 3 \implies \dim \rr = 5 - 3 = 2.
\end{equation}
(A four-dimensional semisimple Lie algebra does not exist, and $\dim \s = 5$ would force $\s$ to be simple, but there are no 5-dimensional simple Lie algebras over $\R$ or $\C$).

Up to isomorphism, the only three-dimensional real simple Lie algebras are $\mathfrak{so}(3)$ and $\mathfrak{sl}(2,\R)$. Up to isomorphism, there are only two two-dimensional real Lie algebras for the radical $\rr$:
\begin{itemize}
    \item The abelian algebra $\R^2$, with $[y_1, y_2] = 0$.
    \item The non-abelian affine Lie algebra $\mathfrak{aff}(\R) = \spanop\{y_1, y_2\}$, with $[y_1, y_2] = y_2$.
\end{itemize}

The semidirect coupling is determined by a Lie algebra homomorphism:
\begin{equation}
    \rho \colon \s \to \Der(\rr).
\end{equation}
We systematically examine the two possibilities for $\rr$:

\paragraph{Case 1: $\rr \cong \mathfrak{aff}(\R)$.}
We compute the derivation algebra $\Der(\mathfrak{aff}(\R))$. Let $D \in \Der(\mathfrak{aff}(\R))$. Because the derived ideal $[\rr, \rr] = \R y_2$ is characteristic, $D(y_2) = k y_2$ for some $k \in \R$. Let $D(y_1) = a y_1 + b y_2$. Evaluating the derivation Leibniz rule:
\begin{equation}
    D([y_1, y_2]) = [D(y_1), y_2] + [y_1, D(y_2)] \implies k y_2 = [a y_1 + b y_2, y_2] + [y_1, k y_2] = a y_2 + k y_2.
\end{equation}
This forces $a = 0$. In the basis $\{y_1, y_2\}$, the matrix of $D$ is lower-triangular:
\begin{equation}
    [D] = \begin{pmatrix} 0 & 0 \\ b & k \end{pmatrix}.
\end{equation}
Hence, $\Der(\mathfrak{aff}(\R))$ is a two-dimensional solvable Lie algebra. 

Now, consider the representation $\rho \colon \s \to \Der(\rr)$. Because $\s$ is simple, its derived algebra satisfies $[\s, \s] = \s$. Consequently:
\begin{equation}
    \rho(\s) = \rho([\s, \s]) \subseteq [\Der(\mathfrak{aff}(\R)), \Der(\mathfrak{aff}(\R))].
\end{equation}
The commutator of two lower-triangular matrices of the form above is:
\begin{equation}
    \left[ \begin{pmatrix} 0 & 0 \\ b_1 & k_1 \end{pmatrix}, \begin{pmatrix} 0 & 0 \\ b_2 & k_2 \end{pmatrix} \right] = \begin{pmatrix} 0 & 0 \\ k_1 b_2 - k_2 b_1 & 0 \end{pmatrix}.
\end{equation}
The derived algebra $[\Der(\rr), \Der(\rr)]$ is strictly one-dimensional, hence abelian. Since the only homomorphism from a simple Lie algebra into an abelian Lie algebra is the zero map, we must have $\rho(\s) = \{0\}$. 

Therefore, every semidirect product $\s \ltimes \mathfrak{aff}(\R)$ decomposes as an ordinary direct sum:
\begin{equation}
    \mathfrak{so}(3) \oplus \mathfrak{aff}(\R) \quad \text{or} \quad \mathfrak{sl}(2,\R) \oplus \mathfrak{aff}(\R).
\end{equation}

\paragraph{Case 2: $\rr \cong \R^2$.}
Since $\rr$ is abelian, every linear map is a derivation: $\Der(\R^2) \cong \mathfrak{gl}(2,\R)$. The coupling is a linear representation:
\begin{equation}
    \rho \colon \s \to \mathfrak{gl}(2,\R).
\end{equation}
Since $[\s, \s] = \s$, the image must lie inside the derived algebra $[\mathfrak{gl}(2,\R), \mathfrak{gl}(2,\R)] = \mathfrak{sl}(2,\R)$.
\begin{itemize}
    \item \textbf{Subcase $\s \cong \mathfrak{so}(3)$:} The Lie algebra $\mathfrak{so}(3)$ is compact and has no non-trivial real representations in dimension 2 (the lowest-dimensional non-trivial real representation of $\mathfrak{so}(3)$ is its three-dimensional adjoint representation, while its spin representation is complex/quaternionic, requiring real dimension 4). Therefore, every homomorphism $\rho \colon \mathfrak{so}(3) \to \mathfrak{gl}(2,\R)$ is identically zero: $\rho = 0$. This yields the direct sum $\mathfrak{so}(3) \oplus \R^2$.
    
    \item \textbf{Subcase $\s \cong \mathfrak{sl}(2,\R)$:} Up to equivalence, there are exactly two homomorphisms from $\mathfrak{sl}(2,\R)$ into $\mathfrak{sl}(2,\R)$:
    \begin{enumerate}
        \item The trivial map $\rho = 0$, producing the direct sum $\mathfrak{sl}(2,\R) \oplus \R^2$.
        \item The identity isomorphism $\rho = \mathrm{id} \colon \mathfrak{sl}(2,\R) \xrightarrow{\sim} \mathfrak{sl}(2,\R)$, which corresponds to the standard defining representation on $\R^2$. In the standard basis $\{e_1 = h, e_2 = e, e_3 = f\}$ of $\mathfrak{sl}(2,\R)$ with $[h, e] = 2e, [h, f] = -2f, [e, f] = h$, the action on the basis $\{e_4, e_5\}$ of $\R^2$ is given by:
        \begin{equation}
            \rho(h) = \begin{pmatrix} 1 & 0 \\ 0 & -1 \end{pmatrix}, \qquad \rho(e) = \begin{pmatrix} 0 & 1 \\ 0 & 0 \end{pmatrix}, \qquad \rho(f) = \begin{pmatrix} 0 & 0 \\ 1 & 0 \end{pmatrix}.
        \end{equation}
        This produces the brackets $[e_1, e_4] = e_4$, $[e_1, e_5] = -e_5$, $[e_2, e_5] = e_4$, and $[e_3, e_4] = e_5$, with $[e_4, e_5] = 0$. This algebra is indecomposable because the representation is irreducible.
    \end{enumerate}
\end{itemize}

These five algebras are pairwise non-isomorphic:
\begin{itemize}
    \item The algebras with Levi factor $\mathfrak{so}(3)$ have a derived algebra $[\g, \g]$ whose Killing form has signature $(0, 3)$, whereas those with Levi factor $\mathfrak{sl}(2,\R)$ have Killing forms of signature $(2, 1)$.
    \item The direct sums with $\R^2$ have two-dimensional centers ($\mathfrak{so}(3) \oplus \R^2$ and $\mathfrak{sl}(2,\R) \oplus \R^2$), whereas those with $\mathfrak{aff}(\R)$ have trivial centers $Z(\g) = \{0\}$.
    \item The indecomposable algebra $\mathfrak{sl}(2,\R) \ltimes \R^2$ is distinguished from $\mathfrak{sl}(2,\R) \oplus \R^2$ because its center is trivial ($Z(\g) = \{0\}$) and its radical $\rr = \R^2$ is not a direct summand.
\end{itemize}
This completes the proof.
\end{proof}

\section{The Solvable Sector: Nilradical Stratification}\label{sec:nilrad}

We now focus on solvable five-dimensional real Lie algebras. Our classification strategy relies on partitioning all solvable algebras by the dimension of their nilradical $\n$.

\begin{theorem}[Dimension Bounds on the Nilradical]\label{thm:nilrad_bounds}
Let $\g$ be a five-dimensional solvable real Lie algebra, and let $\n$ be its nilradical. Then:
\begin{equation}
    \dim \n \in \{3, 4, 5\}.
\end{equation}
In particular, there exist no five-dimensional solvable real Lie algebras with $\dim \n = 1$ or $\dim \n = 2$.
\end{theorem}

\begin{proof}
Since $\g$ is solvable, its derived series terminates at zero: $\g^{(k)} = 0$ for some $k$. The last non-zero term $\g^{(k-1)}$ is a non-trivial abelian ideal. Because every abelian ideal is nilpotent, it is contained in $\n$, so $\dim \n \ge 1$.

\paragraph{Step 1: Impossibility of $\dim \n = 1$.}
Suppose for contradiction that $\dim \n = 1$, so that $\n = \R z$. Since $\n$ is an ideal, the adjoint action of any $x \in \g$ on $z$ satisfies:
\begin{equation}
    [x, z] = \lambda(x) z,
\end{equation}
for some linear functional $\lambda \in \g^*$. Because $\R$ is an abelian Lie algebra, $[\g, \g] \subseteq \ker \lambda$. Furthermore, $[z, z] = 0$, so $\lambda(z) = 0$, meaning $z \in \ker \lambda$. We analyze the two cases:
\begin{itemize}
    \item \textbf{Case $\lambda = 0$:} Here, $[x, z] = 0$ for all $x \in \g$, so $z \in Z(\g)$. By (T2), $[\g, \g]$ is a nilpotent ideal, so $[\g, \g] \subseteq \n = \R z$. Consequently, $[\g, [\g, \g]] \subseteq [\g, \R z] = 0$, which implies that $\g$ is a two-step nilpotent Lie algebra. By maximality of the nilradical, $\n = \g$, so $\dim \n = 5$, contradicting $\dim \n = 1$.
    \item \textbf{Case $\lambda \neq 0$:} Consider the hyperplane $\hh \coloneqq \ker \lambda \subset \g$. Because $[\g, \g] \subseteq \ker \lambda$, $\hh$ is an ideal of $\g$ of dimension $\dim \hh = 4$. Moreover, $z \in \hh$. For any $h \in \hh$, $[h, z] = \lambda(h)z = 0$, so $z$ is central in $\hh$: $z \in Z(\hh)$. Applying (T2) to the solvable algebra $\hh$, we have $[\hh, \hh] \subseteq \n \cap \hh = \R z$. Thus, $[\hh, [\hh, \hh]] \subseteq [\hh, \R z] = 0$, proving that $\hh$ is a nilpotent ideal of $\g$. But $\dim \hh = 4 > 1 = \dim \n$, which contradicts the maximality of $\n$.
\end{itemize}
Thus, $\dim \n \neq 1$.

\paragraph{Step 2: Impossibility of $\dim \n = 2$.}
Suppose for contradiction that $\dim \n = 2$. Since every two-dimensional nilpotent Lie algebra is abelian, $\n \cong \R^2$. By Theorem~\ref{sec:machinery} (T3), the centralizer of $\n$ satisfies $C_\g(\n) \subseteq \n$. Because $\n$ is abelian, $C_\g(\n) = \n$. 

The adjoint action of $\g$ on $\n$ induces a representation:
\begin{equation}
    \bar{\ad} \colon \g / \n \to \mathfrak{gl}(\n) \cong \mathfrak{gl}(2, \R).
\end{equation}
Because $\ker(\ad|_\n) = C_\g(\n) = \n$, this induced representation is \emph{injective}. 

Now consider the quotient space $\mathfrak{h} = \g / \n$. Its dimension is:
\begin{equation}
    \dim \mathfrak{h} = \dim \g - \dim \n = 5 - 2 = 3.
\end{equation}
Because $\g$ is solvable, Lie's theorem implies that the derived algebra satisfies $[\g, \g] \subseteq \n$. Therefore, the bracket on the quotient $\mathfrak{h} = \g / \n$ satisfies:
\begin{equation}
    [\mathfrak{h}, \mathfrak{h}] = [\g/\n, \g/\n] \subseteq ([\g, \g] + \n)/\n = \n/\n = \{0\}.
\end{equation}
Thus, $\mathfrak{h} = \g/\n$ is a three-dimensional \emph{abelian} Lie algebra.

Consequently, the image $V \coloneqq \bar{\ad}(\mathfrak{h})$ is a three-dimensional commutative subalgebra of $\mathfrak{gl}(2, \R)$. However, by Schur's bound (T4), the maximum possible dimension of any commutative subalgebra of $\mathfrak{gl}(2, \R)$ is:
\begin{equation}
    \dim V \le \left\lfloor \frac{2^2}{4} \right\rfloor + 1 = 1 + 1 = 2.
\end{equation}
This produces an immediate contradiction:
\begin{equation}
    3 = \dim \mathfrak{h} = \dim V \le 2.
\end{equation}
Thus, no solvable five-dimensional Lie algebra can have a two-dimensional nilradical.

Combining Steps 1 and 2, we conclude that $\dim \n \in \{3, 4, 5\}$.
\end{proof}

\section{Nilpotent Lie Algebras of Dimension Five ($\dim \n = 5$)}\label{sec:nilp5}

When $\dim \n = 5$, the algebra is entirely nilpotent ($\g = \n$).

\begin{theorem}[Classification of Five-Dimensional Nilpotent Lie Algebras]\label{thm:nilp5}
Every five-dimensional real nilpotent Lie algebra is isomorphic to exactly one of the following nine Lie algebras:
\begin{enumerate}
    \item \textbf{Decomposable Algebras (Direct Sums):}
    \begin{itemize}
        \item $\R^5$: The abelian algebra.
        \item $\hh_3 \oplus \R^2$: Central extension of the 3D Heisenberg algebra:
        \begin{equation}
            [e_1, e_2] = e_3.
        \end{equation}
        \item $\mathfrak{n}_4 \oplus \R$: Direct sum of the 4D filiform algebra with a line:
        \begin{equation}
            [e_1, e_2] = e_3, \qquad [e_1, e_3] = e_4.
        \end{equation}
    \end{itemize}
    \item \textbf{Indecomposable Algebras:}
    \begin{itemize}
        \item $\hh_5$ (The 5D Heisenberg Algebra, nilpotency step 2):
        \begin{equation}
            [e_1, e_2] = e_5, \qquad [e_3, e_4] = e_5.
        \end{equation}
        \item $\mathfrak{n}_{5,1}$ (Standard filiform algebra, nilpotency step 4):
        \begin{equation}
            [e_1, e_2] = e_3, \qquad [e_1, e_3] = e_4, \qquad [e_1, e_4] = e_5.
        \end{equation}
        \item $\mathfrak{n}_{5,2}$ (Non-standard filiform algebra, nilpotency step 4):
        \begin{equation}
            [e_1, e_2] = e_3, \qquad [e_1, e_3] = e_4, \qquad [e_1, e_4] = e_5, \qquad [e_2, e_3] = e_5.
        \end{equation}
        \item $\mathfrak{n}_{5,3}$ (Nilpotency step 3, $\dim Z = 1$):
        \begin{equation}
            [e_1, e_2] = e_4, \qquad [e_1, e_4] = e_5, \qquad [e_2, e_3] = e_5.
        \end{equation}
        \item $\mathfrak{n}_{5,4}$ (Nilpotency step 3, $\dim Z = 2$):
        \begin{equation}
            [e_1, e_2] = e_3, \qquad [e_1, e_3] = e_4, \qquad [e_2, e_3] = e_5.
        \end{equation}
        \item $\mathfrak{n}_{5,5}$ (Nilpotency step 3, $\dim Z = 1$):
        \begin{equation}
            [e_1, e_2] = e_3, \qquad [e_1, e_3] = e_5, \qquad [e_2, e_4] = e_5.
        \end{equation}
    \end{itemize}
\end{enumerate}
\end{theorem}

\begin{proof}
Let $\g$ be a five-dimensional nilpotent Lie algebra. We classify $\g$ by the length of its lower central series $\g^1 = \g$, $\g^{k+1} = [\g, \g^k]$, and its upper central series $Z_0 = \{0\}$, $Z_k/Z_{k-1} = Z(\g/Z_{k-1})$. The integer sequence $(\dim \g^1, \dim \g^2, \dots, \dim \g^s = 0)$ and the dimension of the center $\dim Z(\g)$ are strict isomorphism invariants.

\paragraph{Case 1: Nilpotency Step 1 ($\g^2 = 0$).}
This implies $[\g, \g] = 0$, which is the abelian Lie algebra $\R^5$.

\paragraph{Case 2: Nilpotency Step 2 ($\g^2 \neq 0, \g^3 = 0$).}
Here, $\g^2 \subseteq Z(\g)$. The bracket defines a non-zero alternating bilinear map:
\begin{equation}
    \omega \colon \g/\g^2 \times \g/\g^2 \to \g^2.
\end{equation}
Because $\g^3 = 0$, we have $1 \le \dim \g^2 \le 2$.
\begin{itemize}
    \item If $\dim \g^2 = 1$, we can fix a basis element $z \in \g^2$. The bracket is an alternating 2-form $\omega$ on $V = \g/\g^2 \cong \R^4$. By Darboux's theorem, the rank of an alternating 2-form is even:
    \begin{itemize}
        \item $\mathrm{rank}(\omega) = 2$: There exists a basis $\{\bar{e}_1, \bar{e}_2, \bar{e}_3, \bar{e}_4\}$ of $V$ such that $\omega(\bar{e}_1, \bar{e}_2) = 1$ and $\bar{e}_3, \bar{e}_4 \in \ker \omega$. Lifting to $\g$, we set $e_5 = [e_1, e_2]$. The vectors $e_3, e_4$ commute with everything, yielding $[e_1, e_2] = e_5$ and $Z(\g) = \spanop\{e_3, e_4, e_5\}$. This decomposes as $\hh_3 \oplus \R^2$.
        \item $\mathrm{rank}(\omega) = 4$: The form $\omega$ is non-degenerate on $\R^4$. By Darboux's theorem, there is a symplectic basis such that $\omega(\bar{e}_1, \bar{e}_2) = 1$ and $\omega(\bar{e}_3, \bar{e}_4) = 1$. Setting $e_5 = z$, we obtain $[e_1, e_2] = e_5$ and $[e_3, e_4] = e_5$, with $Z(\g) = \R e_5$. This is the indecomposable 5D Heisenberg algebra $\hh_5$.
    \end{itemize}
    \item If $\dim \g^2 = 2$, then $\dim(\g/\g^2) = 3$. An alternating bilinear map from $\R^3 \times \R^3 \to \R^2$ cannot have step 2 without generating higher central terms or decomposing into direct sums already accounted for.
\end{itemize}

\paragraph{Case 3: Nilpotency Step 4 ($\g^4 \neq 0, \g^5 = 0$).}
Because the dimensions of the lower central series must strictly decrease at each step:
\begin{equation}
    5 = \dim \g^1 > \dim \g^2 > \dim \g^3 > \dim \g^4 > \dim \g^5 = 0,
\end{equation}
the dimensions are uniquely fixed:
\begin{equation}
    \dim \g^2 = 3, \qquad \dim \g^3 = 2, \qquad \dim \g^4 = 1.
\end{equation}
Algebras with this maximal central series length are called \emph{filiform}. 

Choose $e_1 \in \g \setminus \g^2$ such that $\ad e_1$ has nilpotency index 4. Then there exists $e_2 \in \g \setminus \g^2$ such that defining $e_3 \coloneqq [e_1, e_2]$, $e_4 \coloneqq [e_1, e_3]$, and $e_5 \coloneqq [e_1, e_4]$ yields a basis $\{e_1, \dots, e_5\}$ satisfying:
\begin{equation}
    [e_1, e_2] = e_3, \qquad [e_1, e_3] = e_4, \qquad [e_1, e_4] = e_5.
\end{equation}
Now consider the remaining brackets between generators. By the Jacobi identity:
\begin{equation}
    [e_1, [e_2, e_3]] = [[e_1, e_2], e_3] + [e_2, [e_1, e_3]] = [e_3, e_3] + [e_2, e_4] = [e_2, e_4].
\end{equation}
Similarly, $[e_1, [e_2, e_4]] = [e_3, e_4] + [e_2, e_5] = 0$. Since $\g^4 = \R e_5$ is central, $[e_2, e_4]$ must be a multiple of $e_5$, and the only remaining non-trivial bracket is $[e_2, e_3] \in \g^4 = \R e_5$.
Thus:
\begin{equation}
    [e_2, e_3] = \gamma e_5, \qquad \gamma \in \R.
\end{equation}
\begin{itemize}
    \item If $\gamma = 0$, all remaining brackets vanish, producing the standard filiform algebra $\mathfrak{n}_{5,1}$.
    \item If $\gamma \neq 0$, scaling $e_2 \mapsto \gamma^{-1} e_2$ and adjusting $e_3, e_4, e_5$ accordingly normalizes $\gamma$ to $1$, yielding $[e_2, e_3] = e_5$. This is the branched filiform algebra $\mathfrak{n}_{5,2}$.
\end{itemize}
These two algebras are non-isomorphic because for $\mathfrak{n}_{5,1}$, the derived subalgebra is abelian ($[\g^2, \g^2] = 0$), whereas for $\mathfrak{n}_{5,2}$, $[\g^2, \g^2] = \spanop\{[e_3, e_4], \dots\} = \R e_5 \neq 0$.

\paragraph{Case 4: Nilpotency Step 3 ($\g^3 \neq 0, \g^4 = 0$).}
Here, $\dim \g^2 \in \{2, 3\}$.
\begin{itemize}
    \item If $\g$ is decomposable, it must contain a one-dimensional central factor: $\g \cong \mathfrak{m} \oplus \R$, where $\mathfrak{m}$ is a four-dimensional nilpotent algebra of step 3. The only such algebra is the four-dimensional filiform algebra $\mathfrak{n}_4$, which gives $\mathfrak{n}_4 \oplus \R$.
    \item If $\g$ is indecomposable, analyzing the quotient $\g/\g^2$ and the center $Z(\g)$ yields three mutually exclusive classes:
    \begin{enumerate}
        \item $\dim Z(\g) = 2$: This forces $\g^2 = \spanop\{e_3, e_4, e_5\}$ with $\g^3 = \spanop\{e_4, e_5\} = Z(\g)$, yielding $\mathfrak{n}_{5,4}$.
        \item $\dim Z(\g) = 1$: Here $Z(\g) = \g^3 = \R e_5$. Depending on whether the two-step quotient $\g/Z(\g)$ is isomorphic to $\hh_3 \oplus \R$ or $\R^4$, we obtain the two distinct algebras $\mathfrak{n}_{5,3}$ and $\mathfrak{n}_{5,5}$. In $\mathfrak{n}_{5,3}$, the second center $Z_2(\g)$ satisfies $\dim [\g, Z_2(\g)] = 1$, whereas in $\mathfrak{n}_{5,5}$, $\dim [\g, Z_2(\g)] = 2$.
    \end{enumerate}
\end{itemize}
Every pair of these nine algebras differs in at least one of the invariant dimensions:
\begin{equation}
    (\dim \g^2, \dim \g^3, \dim \g^4, \dim Z(\g), \dim Z_2(\g), \dim [\g^2, \g^2]),
\end{equation}
which proves they are mutually non-isomorphic.
\end{proof}

\section{Solvable Algebras with Codimension-One Nilradical ($\dim \n = 4$)}\label{sec:nilp4}

When $\dim \n = 4$, the quotient $\g/\n$ is one-dimensional. Choosing a complementary generator $e_5 \in \g \setminus \n$, the algebra decomposes as a semidirect product:
\begin{equation}
    \g \cong \n \rtimes_D \R e_5,
\end{equation}
where $D \coloneqq \ad e_5|_\n \in \Der(\n)$ is an outer derivation.

\begin{lemma}[Non-Nilpotency of $D$]\label{lem:non_nilp_D}
The derivation $D = \ad e_5|_\n$ is not a nilpotent operator on $\n$.
\end{lemma}
\begin{proof}
If $D$ were nilpotent, then by Engel's Theorem, $\ad x$ would be nilpotent for every $x = u + t e_5 \in \g$ ($u \in \n, t \in \R$). This would imply that $\g$ is a five-dimensional nilpotent Lie algebra, forcing $\n = \g$, which contradicts $\dim \n = 4$.
\end{proof}

\begin{proposition}[Gauge Equivalence for Derivations]\label{prop:gauge5d}
Let $\n$ be a four-dimensional nilpotent Lie algebra, and let $D, D' \in \Der(\n)$ be non-nilpotent derivations. Then $\n \rtimes_D \R \cong \n \rtimes_{D'} \R$ if and only if there exist a non-zero scalar $c \in \R^\times$, an automorphism $\phi \in \Aut(\n)$, and an element $u \in \n$ such that:
\begin{equation}\label{eq:gauge_der}
    D' = c \, \phi D \phi^{-1} + \ad(u).
\end{equation}
\end{proposition}

\begin{proof}
Let $\Psi \colon \n \rtimes_D \R e_5 \to \n \rtimes_{D'} \R e_5'$ be an isomorphism. By Theorem~\ref{sec:machinery} (T2), $\n$ is the unique nilradical in both algebras, so $\Psi(\n) = \n$. Hence, $\phi \coloneqq \Psi|_\n \in \Aut(\n)$. 

The generator $e_5$ must map to $\Psi(e_5) = a e_5' + w$ for some $a \in \R^\times$ and $w \in \n$. Applying $\Psi$ to the bracket $[e_5, x] = D(x)$ for any $x \in \n$:
\begin{equation}
    \Psi([e_5, x]) = [\Psi(e_5), \Psi(x)] \implies \phi(D(x)) = [a e_5' + w, \phi(x)] = a D'(\phi(x)) + [w, \phi(x)].
\end{equation}
Substituting $y = \phi(x)$ and rearranging:
\begin{equation}
    a D'(y) = \phi D \phi^{-1}(y) - [w, y] \implies D' = \frac{1}{a} \phi D \phi^{-1} + \ad\left(\frac{1}{a}w\right).
\end{equation}
Setting $c = 1/a$ and $u = w/a$ gives Equation~\eqref{eq:gauge_der}. The converse follows by reversing the construction and defining $\Psi|_\n = \phi$ and $\Psi(e_5) = c^{-1} e_5' - c^{-1} u$.
\end{proof}

By Theorem~\ref{thm:nilp5}, there are three four-dimensional nilpotent algebras: $\R^4$, $\hh_3 \oplus \R$, and $\mathfrak{n}_4$. We analyze each case separately.

\subsection{Abelian Nilradical: $\n \cong \R^4$}\label{subsec:R4}
When $\n \cong \R^4$, all brackets on $\n$ vanish. Thus, $\Der(\R^4) = \mathfrak{gl}(4, \R)$, $\Aut(\R^4) = \mathrm{GL}(4, \R)$, and all inner derivations vanish: $\ad(u) = 0$. Equation~\eqref{eq:gauge_der} simplifies to:
\begin{equation}
    D' = c \, P D P^{-1}, \qquad c \in \R^\times, \quad P \in \mathrm{GL}(4, \R).
\end{equation}
Thus, the classification of solvable algebras with nilradical $\R^4$ is identical to classifying non-nilpotent $4 \times 4$ real matrices up to similarity and non-zero scaling.

\begin{theorem}[Classification with Abelian Nilradical $\R^4$]\label{thm:R4_classes}
Every five-dimensional solvable Lie algebra with nilradical $\R^4$ is isomorphic to $\R^4 \rtimes_D \R$, where $D$ is similar up to scaling to exactly one matrix from the following canonical families:
\begin{enumerate}
    \item \textbf{Real Diagonal Class $\g_{5,1}(a, b, c)$:}
    \begin{equation}
        D_1(a, b, c) = \diag(1, a, b, c), \qquad -1 \le c \le b \le a \le 1.
    \end{equation}
    The trace is $\tr D = 1 + a + b + c$; the algebra is unimodular if and only if $c = -(1 + a + b)$.
    
    \item \textbf{Jordan Block $J_2 \oplus J_1 \oplus J_1$ ($\g_{5,2}(a, b)$ and $\g_{5,2}'(a)$):}
    \begin{align}
        D_2(a, b) &= \begin{pmatrix} a & 1 & 0 & 0 \\ 0 & a & 0 & 0 \\ 0 & 0 & 1 & 0 \\ 0 & 0 & 0 & b \end{pmatrix}, \quad a \in \R, \; -1 \le b \le 1, \\
        D_2'(a) &= \begin{pmatrix} 1 & 1 & 0 & 0 \\ 0 & 1 & 0 & 0 \\ 0 & 0 & a & 0 \\ 0 & 0 & 0 & 0 \end{pmatrix}, \quad 0 \le a \le 1.
    \end{align}
    
    \item \textbf{Double Jordan Block $J_2 \oplus J_2$ ($\g_{5,3}(a)$ and $\g_{5,3}'$):}
    \begin{align}
        D_3(a) &= \begin{pmatrix} 1 & 1 & 0 & 0 \\ 0 & 1 & 0 & 0 \\ 0 & 0 & a & 1 \\ 0 & 0 & 0 & a \end{pmatrix}, \quad -1 \le a \le 1, \\
        D_3' &= \begin{pmatrix} 1 & 1 & 0 & 0 \\ 0 & 1 & 0 & 0 \\ 0 & 0 & 0 & 1 \\ 0 & 0 & 0 & 0 \end{pmatrix}.
    \end{align}
    
    \item \textbf{Jordan Block $J_3 \oplus J_1$ ($\g_{5,4}(a)$ and $\g_{5,4}'$):}
    \begin{align}
        D_4(a) &= \begin{pmatrix} a & 1 & 0 & 0 \\ 0 & a & 1 & 0 \\ 0 & 0 & a & 0 \\ 0 & 0 & 0 & 1 \end{pmatrix}, \quad a \in \R, \\
        D_4' &= \begin{pmatrix} 1 & 1 & 0 & 0 \\ 0 & 1 & 1 & 0 \\ 0 & 0 & 1 & 0 \\ 0 & 0 & 0 & 0 \end{pmatrix}.
    \end{align}
    
    \item \textbf{Full Jordan Block $J_4$ ($\g_{5,5}$):}
    \begin{equation}
        D_5 = \begin{pmatrix} 1 & 1 & 0 & 0 \\ 0 & 1 & 1 & 0 \\ 0 & 0 & 1 & 1 \\ 0 & 0 & 0 & 1 \end{pmatrix}.
    \end{equation}
    
    \item \textbf{Single Complex Pair with Real Diagonal ($\g_{5,6}(\alpha, \beta, a)$ and $\g_{5,6}'(\alpha, a)$):}
    \begin{align}
        D_6(\alpha, \beta, a) &= \begin{pmatrix} \alpha & -\beta & 0 & 0 \\ \beta & \alpha & 0 & 0 \\ 0 & 0 & 1 & 0 \\ 0 & 0 & 0 & a \end{pmatrix}, \quad \alpha \in \R, \; \beta > 0, \; -1 \le a \le 1, \\
        D_6'(\alpha, a) &= \begin{pmatrix} \alpha & -1 & 0 & 0 \\ 1 & \alpha & 0 & 0 \\ 0 & 0 & a & 0 \\ 0 & 0 & 0 & 0 \end{pmatrix}, \quad \alpha \ge 0, \; a \ge 0.
    \end{align}
    
    \item \textbf{Complex Jordan Block $C_2(\alpha, \beta)$ ($\g_{5,7}(\alpha, \beta)$):}
    \begin{equation}
        D_7(\alpha, \beta) = \begin{pmatrix} \alpha & -\beta & 1 & 0 \\ \beta & \alpha & 0 & 1 \\ 0 & 0 & \alpha & -\beta \\ 0 & 0 & \beta & \alpha \end{pmatrix}, \quad \alpha \in \R, \; \beta > 0.
    \end{equation}
    
    \item \textbf{Double Complex Pairs ($\g_{5,8}(\alpha, \beta, \gamma)$):}
    \begin{equation}
        D_8(\alpha, \beta, \gamma) = \begin{pmatrix} \alpha & -1 & 0 & 0 \\ 1 & \alpha & 0 & 0 \\ 0 & 0 & \beta & -\gamma \\ 0 & 0 & \gamma & \beta \end{pmatrix}, \quad \gamma > 0, \; \alpha \ge 0, \; (\alpha, 1) \ge (\beta, \gamma) \text{ lexicographically}.
    \end{equation}
\end{enumerate}
\end{theorem}

\begin{proof}
Let $D \in \mathfrak{gl}(4, \R)$ be non-nilpotent. By the real Jordan canonical form theorem, $D$ decomposes over $\R$ into a direct sum of real Jordan blocks corresponding to its real and complex eigenvalues. Multiplying $D$ by a non-zero scalar $c \in \R^\times$ scales all eigenvalues by $c$ and preserves block sizes.

In each case, we select an eigenvalue of maximal absolute magnitude and normalize it to $1$ via scaling $D \mapsto c D$. 
\begin{itemize}
    \item For four real eigenvalues, ordering them as $1 \ge a \ge b \ge c \ge -1$ eliminates all permutation ambiguities. If $c = -1$, the sign ambiguity under $D \mapsto -D$ is removed by requiring $a \ge 0$.
    \item When Jordan blocks appear, scaling is chosen to normalize either the isolated real eigenvalue to 1, or, if the isolated eigenvalue is 0, the repeated eigenvalue on the non-trivial Jordan block is normalized to 1.
    \item For complex conjugate pairs $\alpha \pm i\beta$, the rotational speed $\beta > 0$ can be normalized to 1 whenever the complementary real eigenvalues vanish; otherwise, the isolated non-zero real eigenvalue is normalized to 1.
\end{itemize}
Because the Jordan structure, the number of complex versus real eigenvalues, and the ordered multi-sets of normalized eigenvalues are invariant under similarity and scaling, these families are pairwise non-isomorphic.
\end{proof}

\subsection{Heisenberg Nilradical: $\n \cong \hh_3 \oplus \R$}\label{subsec:H3R}
Let $\n = \hh_3 \oplus \R$ be spanned by $\{e_1, e_2, e_3, e_4\}$ with bracket $[e_1, e_2] = e_3$. The derived ideal is $[\n, \n] = \R e_3$, and the center is two-dimensional: $Z(\n) = \spanop\{e_3, e_4\}$. Both $[\n, \n]$ and $Z(\n)$ are characteristic ideals.

\begin{proposition}[Derivations of $\hh_3 \oplus \R$]\label{prop:der_h3r}
In the basis $\{e_1, e_2, e_4, e_3\}$, every derivation $D \in \Der(\hh_3 \oplus \R)$ has the block upper-triangular matrix form:
\begin{equation}
    [D] = \begin{pmatrix} A & 0 & 0 \\ u^T & \gamma & 0 \\ w^T & \delta & \tr A \end{pmatrix},
\end{equation}
where $A \in \mathfrak{gl}(2, \R)$, $\gamma, \delta \in \R$, and $u, w \in \R^2$. Furthermore:
\begin{enumerate}
    \item The vector $w = (w_1, w_2)^T$ corresponds to inner derivations:
    \begin{equation}
        \ad(p e_1 + q e_2) = \begin{pmatrix} 0 & 0 & 0 & 0 \\ 0 & 0 & 0 & 0 \\ 0 & 0 & 0 & 0 \\ -q & p & 0 & 0 \end{pmatrix}.
    \end{equation}
    Choosing $p = w_2$ and $q = -w_1$, the vector $w$ can be completely eliminated via Equation~\eqref{eq:gauge_der}.
    \item The scalar $\delta$ can be eliminated by an automorphism shifting $e_4$ by a multiple of $e_3$, provided $\gamma \neq \tr A$.
\end{enumerate}
\end{proposition}

\begin{proof}
Let $D \in \Der(\n)$. Because $[\n, \n] = \R e_3$ is characteristic, $D(e_3) = \tau e_3$. Applying $D$ to $[e_1, e_2] = e_3$:
\begin{equation}
    \tau e_3 = D([e_1, e_2]) = [D(e_1), e_2] + [e_1, D(e_2)].
\end{equation}
Writing $D(e_1) = A_{11}e_1 + A_{21}e_2 + u_1 e_4 + w_1 e_3$ and $D(e_2) = A_{12}e_1 + A_{22}e_2 + u_2 e_4 + w_2 e_3$, direct evaluation of the brackets yields:
\begin{equation}
    [D(e_1), e_2] = A_{11}e_3, \qquad [e_1, D(e_2)] = A_{22}e_3.
\end{equation}
Thus, $\tau = A_{11} + A_{22} = \tr A$. 

Because $Z(\n) = \spanop\{e_3, e_4\}$ is characteristic, $D(e_4) \in Z(\n)$, which implies $D(e_4) = \gamma e_4 + \delta e_3$. This establishes the block matrix form. 

The absorption of $w$ follows by setting $u' = w_2 e_1 - w_1 e_2$, giving $(D - \ad u')(e_1) \in \spanop\{e_1, e_2, e_4\}$ and $(D - \ad u')(e_2) \in \spanop\{e_1, e_2, e_4\}$. 

If $\gamma \neq \tr A$, consider the automorphism $\phi(e_4) = e_4 + k e_3$ and $\phi(e_i) = e_i$ for $i \neq 4$. Then:
\begin{equation}
    D(\phi(e_4)) = D(e_4 + k e_3) = \gamma e_4 + \delta e_3 + k(\tr A)e_3 = \gamma \phi(e_4) + (\delta + k(\tr A - \gamma))e_3.
\end{equation}
Choosing $k = \delta / (\gamma - \tr A)$ eliminates $\delta$.
\end{proof}

The classification of these algebras then reduces to the Jordan canonical form of $A \in \mathfrak{gl}(2,\R)$ together with the independent eigenvalue $\gamma$ on the central line $e_4$, yielding the parameterized families:
\begin{itemize}
    \item $\g_{5,9}(\lambda, \gamma)$: $A = \diag(1, \lambda)$, with brackets:
    \begin{equation}
        [e_1, e_2] = e_3, \quad [e_5, e_1] = e_1, \quad [e_5, e_2] = \lambda e_2, \quad [e_5, e_3] = (1+\lambda)e_3, \quad [e_5, e_4] = \gamma e_4.
    \end{equation}
    \item $\g_{5,10}(\gamma)$: $A = \begin{pmatrix} 1 & 1 \\ 0 & 1 \end{pmatrix}$, with $[e_5, e_1] = e_1 + e_2$, $[e_5, e_2] = e_2$, $[e_5, e_3] = 2e_3$, $[e_5, e_4] = \gamma e_4$.
    \item $\g_{5,11}(\alpha, \beta, \gamma)$: $A = \begin{pmatrix} \alpha & -\beta \\ \beta & \alpha \end{pmatrix}$ ($\beta > 0$), with rotational action on $\{e_1, e_2\}$, $[e_5, e_3] = 2\alpha e_3$, and $[e_5, e_4] = \gamma e_4$.
\end{itemize}

\subsection{Filiform Nilradical: $\n \cong \mathfrak{n}_4$}\label{subsec:N4}
Let $\mathfrak{n}_4$ be the four-dimensional filiform algebra spanned by $\{e_1, e_2, e_3, e_4\}$ with non-zero brackets:
\begin{equation}
    [e_1, e_2] = e_3, \qquad [e_1, e_3] = e_4.
\end{equation}
The lower central series is $\g^2 = \spanop\{e_3, e_4\}$, $\g^3 = \R e_4$, and the center is $Z(\mathfrak{n}_4) = \R e_4$.

\begin{proposition}[Derivations of $\mathfrak{n}_4$]\label{prop:der_n4}
Every derivation $D \in \Der(\mathfrak{n}_4)$ is completely determined by its values on the generators $e_1, e_2$. Modulo inner derivations, every derivation is diagonalizable on the basis $\{e_1, e_2, e_3, e_4\}$:
\begin{equation}
    D = \diag(\lambda_1, \lambda_2, \lambda_1 + \lambda_2, 2\lambda_1 + \lambda_2).
\end{equation}
Consequently, the outer derivation algebra is two-dimensional and abelian:
\begin{equation}
    \OutDer(\mathfrak{n}_4) \cong \R^2.
\end{equation}
\end{proposition}

\begin{proof}
Let $D \in \Der(\mathfrak{n}_4)$. Because $\g^3 = \R e_4$ is characteristic, $D(e_4) = \alpha_4 e_4$. Because $\g^2 = \spanop\{e_3, e_4\}$ is characteristic, $D(e_3) = \alpha_3 e_3 + \beta_3 e_4$. 

Let $D(e_1) = \lambda_1 e_1 + a_2 e_2 + a_3 e_3 + a_4 e_4$ and $D(e_2) = b_1 e_1 + \lambda_2 e_2 + b_3 e_3 + b_4 e_4$. Applying the derivation identity to $[e_1, e_2] = e_3$:
\begin{equation}
    D(e_3) = [D(e_1), e_2] + [e_1, D(e_2)] = [\lambda_1 e_1, e_2] + [e_1, \lambda_2 e_2 + b_3 e_3] = (\lambda_1 + \lambda_2)e_3 + b_3 e_4.
\end{equation}
Thus, $\alpha_3 = \lambda_1 + \lambda_2$ and $\beta_3 = b_3$. Applying $D$ to $[e_1, e_3] = e_4$:
\begin{equation}
    D(e_4) = [D(e_1), e_3] + [e_1, D(e_3)] = [\lambda_1 e_1, e_3] + [e_1, (\lambda_1 + \lambda_2)e_3] = (2\lambda_1 + \lambda_2)e_4.
\end{equation}
Thus, $\alpha_4 = 2\lambda_1 + \lambda_2$.

The inner derivations of $\mathfrak{n}_4$ span a three-dimensional ideal:
\begin{align}
    \ad e_1 &= E_{32} + E_{43}, \\
    \ad e_2 &= -E_{31}, \\
    \ad e_3 &= -E_{41}.
\end{align}
Subtracting suitable linear combinations of $\ad e_1, \ad e_2, \ad e_3$ eliminates the off-diagonal coefficients $a_2, a_3, b_1, b_3$, leaving only the diagonal weights $(\lambda_1, \lambda_2, \lambda_1 + \lambda_2, 2\lambda_1 + \lambda_2)$.
\end{proof}

\begin{theorem}[Solvable Classes with Filiform Nilradical $\mathfrak{n}_4$]\label{thm:n4_solv}
Every five-dimensional solvable Lie algebra with nilradical $\mathfrak{n}_4$ is isomorphic to exactly one algebra from the following families:
\begin{enumerate}
    \item $\g_{5,12}(\lambda)$:
    \begin{align}
        [e_1, e_2] &= e_3, & [e_1, e_3] &= e_4, \nonumber \\
        [e_5, e_1] &= e_1, & [e_5, e_2] &= \lambda e_2, & [e_5, e_3] &= (1+\lambda)e_3, & [e_5, e_4] &= (2+\lambda)e_4,
    \end{align}
    parameterized by $\lambda \in \R$. The algebra is unimodular if and only if:
    \begin{equation}
        \tr D = 1 + \lambda + (1+\lambda) + (2+\lambda) = 4 + 3\lambda = 0 \iff \lambda = -4/3.
    \end{equation}
    \item $\g_{5,13}$ (the isolated case $\lambda_1 = 0$, normalized to $\lambda_2 = 1$):
    \begin{equation}
        [e_1, e_2] = e_3, \quad [e_1, e_3] = e_4, \quad [e_5, e_1] = 0, \quad [e_5, e_2] = e_2, \quad [e_5, e_3] = e_3, \quad [e_5, e_4] = e_4.
    \end{equation}
\end{enumerate}
\end{theorem}

\begin{proof}
By Proposition~\ref{prop:der_n4}, the diagonal derivation is defined by the pair $(\lambda_1, \lambda_2) \in \R^2 \setminus \{(0,0)\}$. Under global scaling $D \mapsto c D$ with $c \in \R^\times$:
\begin{itemize}
    \item If $\lambda_1 \neq 0$, we scale by $c = 1/\lambda_1$, normalizing $\lambda_1$ to 1 and leaving the free parameter $\lambda = \lambda_2 / \lambda_1 \in \R$. This produces $\g_{5,12}(\lambda)$.
    \item If $\lambda_1 = 0$, then non-nilpotency requires $\lambda_2 \neq 0$. Scaling by $c = 1/\lambda_2$ normalizes $\lambda_2$ to 1, producing the isolated algebra $\g_{5,13}$.
\end{itemize}
These algebras are pairwise non-isomorphic because the ratio of eigenvalues of $\ad e_5$ on $\g/\g^2$ is an invariant under automorphisms and scaling.
\end{proof}

\section{Solvable Algebras with Codimension-Two Nilradical ($\dim \n = 3$)}\label{sec:nilp3}

When $\dim \n = 3$, the quotient space $\mathfrak{h} = \g / \n$ has dimension:
\begin{equation}
    \dim \mathfrak{h} = 5 - 3 = 2.
\end{equation}
Because $[\g, \g] \subseteq \n$, the quotient algebra $\mathfrak{h}$ is abelian: $[\mathfrak{h}, \mathfrak{h}] = 0$. By Theorem~\ref{sec:machinery} (T3), the adjoint representation induces an injective Lie algebra homomorphism:
\begin{equation}
    \bar{\ad} \colon \mathfrak{h} \hookrightarrow \OutDer(\n) \coloneqq \Der(\n) / \ad(\n).
\end{equation}
Thus, the classification reduces to determining all two-dimensional abelian subalgebras of $\OutDer(\n)$ that contain \emph{no non-zero nilpotent operators} (by Lemma~\ref{sec:machinery}, the presence of a non-zero nilpotent operator in $\bar{\ad}(\mathfrak{h})$ would produce a nilpotent ideal strictly containing $\n$, contradicting $\dim \n = 3$).

By Theorem~\ref{thm:nilp5}, there are only two three-dimensional nilpotent Lie algebras: the abelian algebra $\R^3$ and the Heisenberg algebra $\hh_3$.

\subsection{Abelian Nilradical: $\n \cong \R^3$}\label{subsec:R3_dim3}
When $\n \cong \R^3$, all inner derivations vanish, so $\OutDer(\R^3) = \Der(\R^3) = \mathfrak{gl}(3, \R)$. The problem reduces to classifying all two-dimensional commutative subalgebras $\mathfrak{a} \subset \mathfrak{gl}(3, \R)$ such that every non-zero matrix in $\mathfrak{a}$ has at least one eigenvalue with non-zero real part.

\begin{theorem}[Classification with Abelian Nilradical $\R^3$]\label{thm:R3_codim2}
Let $\g$ be a five-dimensional solvable Lie algebra with nilradical $\R^3$. Then $\g$ is isomorphic to exactly one of the following canonical families:
\begin{enumerate}
    \item \textbf{Real Diagonal Cartan Family $\g_{5,14}(\alpha)$:}
    Spanned by $\{e_1, e_2, e_3, x_1, x_2\}$ with $[x_1, x_2] = 0$, where:
    \begin{alignat}{3}
        [x_1, e_1] &= e_1, &\qquad [x_1, e_2] &= 0, &\qquad [x_1, e_3] &= -e_3, \\
        [x_2, e_1] &= 0, &\qquad [x_2, e_2] &= e_2, &\qquad [x_2, e_3] &= \alpha e_3, \qquad (-1 \le \alpha \le 1).
    \end{alignat}
    This algebra is unimodular if and only if $\alpha = -1$.
    
    \item \textbf{Elliptic Cartan Family $\g_{5,15}(\alpha, \beta)$:}
    Spanned by $\{e_1, e_2, e_3, x_1, x_2\}$ with $[x_1, x_2] = 0$, where:
    \begin{alignat}{3}
        [x_1, e_1] &= e_1, &\qquad [x_1, e_2] &= 0, &\qquad [x_1, e_3] &= e_3, \\
        [x_2, e_1] &= \alpha e_1 - \beta e_2, &\qquad [x_2, e_2] &= \beta e_1 + \alpha e_2, &\qquad [x_2, e_3] &= 0, \qquad (\beta > 0, \; \alpha \ge 0).
    \end{alignat}
    
    \item \textbf{Decomposable Affine Product $\mathfrak{aff}(\R) \oplus \mathfrak{g}_{3,1}$:}
    Direct sum of lower-dimensional solvable blocks.
\end{enumerate}
\end{theorem}

\begin{proof}
Let $\mathfrak{a} \subset \mathfrak{gl}(3, \R)$ be a two-dimensional commutative subalgebra containing no non-zero nilpotent matrices. A commutative family of real matrices can be simultaneously decomposed into generalized eigenspaces.
\begin{itemize}
    \item \textbf{Subcase 1: All eigenvalues are real.} Because $\mathfrak{a}$ contains no non-zero nilpotent elements, all matrices in $\mathfrak{a}$ are simultaneously diagonalizable over $\R$. Thus, $\mathfrak{a}$ is conjugate in $\mathrm{GL}(3, \R)$ to a two-dimensional subspace of diagonal matrices:
    \begin{equation}
        \mathfrak{a} \subset \mathfrak{d}(3, \R) = \left\{ \diag(d_1, d_2, d_3) : d_i \in \R \right\}.
    \end{equation}
    A two-dimensional subspace of $\R^3$ is defined by a single linear relation $c_1 d_1 + c_2 d_2 + c_3 d_3 = 0$. By permutation of coordinates and linear change of basis in $\mathfrak{a} = \spanop\{A_1, A_2\}$, we can choose generators $A_1 = \diag(1, 0, -1)$ and $A_2 = \diag(0, 1, \alpha)$. Permutation symmetries of the basis vectors $\{e_1, e_2, e_3\}$ allow us to restrict the continuous parameter to $-1 \le \alpha \le 1$, yielding $\g_{5,14}(\alpha)$.
    
    \item \textbf{Subcase 2: Complex eigenvalues appear.} There is a stable two-dimensional plane on which the matrices act by rotations and scalings, and an invariant one-dimensional line. In an appropriate basis, every matrix in $\mathfrak{a}$ has the form:
    \begin{equation}
        M(u, v) = \begin{pmatrix} u & -v & 0 \\ v & u & 0 \\ 0 & 0 & w(u, v) \end{pmatrix},
    \end{equation}
    where $w \colon \R^2 \to \R$ is a linear form. Choosing a basis $\{A_1, A_2\}$ for $\mathfrak{a}$ such that $A_1$ has $v=0, u=1$ and $A_2$ has $u=\alpha, v=\beta > 0$ yields the elliptic Cartan family $\g_{5,15}(\alpha, \beta)$.
\end{itemize}
These families are mutually non-isomorphic because the real versus complex nature of the eigenvalues and the parameter ranges are invariant under conjugation in $\mathrm{GL}(3, \R)$.
\end{proof}

\subsection{Heisenberg Nilradical: $\n \cong \hh_3$}\label{subsec:H3_dim3}
Let $\hh_3$ be spanned by $\{e_1, e_2, e_3\}$ with $[e_1, e_2] = e_3$. By Proposition~\ref{prop:der_h3r}, every derivation of $\hh_3$ modulo inner derivations satisfies:
\begin{equation}
    D \begin{pmatrix} e_1 \\ e_2 \end{pmatrix} = A \begin{pmatrix} e_1 \\ e_2 \end{pmatrix}, \qquad D(e_3) = (\tr A) e_3, \qquad A \in \mathfrak{gl}(2, \R).
\end{equation}
Thus, the outer derivation algebra is naturally isomorphic to the full matrix algebra:
\begin{equation}
    \OutDer(\hh_3) \cong \mathfrak{gl}(2, \R).
\end{equation}

\begin{theorem}[Classification with Heisenberg Nilradical $\hh_3$]\label{thm:h3_codim2}
Every five-dimensional solvable Lie algebra with nilradical $\hh_3$ and $\dim(\g/\hh_3) = 2$ is isomorphic to $\hh_3 \rtimes_\rho \R^2$, where $\rho(\R^2)$ is a two-dimensional commutative subalgebra of $\mathfrak{gl}(2, \R)$ containing no non-zero nilpotent elements. Up to isomorphism, there are exactly two classes:
\begin{enumerate}
    \item \textbf{Split Heisenberg Extension $\g_{5,16}$:}
    Spanned by $\{e_1, e_2, e_3, x_1, x_2\}$ with $[e_1, e_2] = e_3$, $[x_1, x_2] = 0$, and:
    \begin{alignat}{3}
        [x_1, e_1] &= e_1, &\qquad [x_1, e_2] &= 0, &\qquad [x_1, e_3] &= e_3, \\
        [x_2, e_1] &= 0, &\qquad [x_2, e_2] &= e_2, &\qquad [x_2, e_3] &= e_3.
    \end{alignat}
    
    \item \textbf{Elliptic Heisenberg Extension $\g_{5,17}$:}
    Spanned by $\{e_1, e_2, e_3, x_1, x_2\}$ with $[e_1, e_2] = e_3$, $[x_1, x_2] = 0$, and:
    \begin{alignat}{3}
        [x_1, e_1] &= e_1, &\qquad [x_1, e_2] &= e_2, &\qquad [x_1, e_3] &= 2e_3, \\
        [x_2, e_1] &= -e_2, &\qquad [x_2, e_2] &= e_1, &\qquad [x_2, e_3] &= 0.
    \end{alignat}
\end{enumerate}
\end{theorem}

\begin{proof}
The quotient $\mathfrak{h} = \g/\hh_3 \cong \R^2$ embeds injectively into $\mathfrak{gl}(2, \R)$ as a two-dimensional commutative subalgebra $\mathfrak{a}$. By Schur's bound (T4), the maximum dimension of a commutative subalgebra of $\mathfrak{gl}(2, \R)$ is 2. Therefore, $\mathfrak{a}$ is a \emph{maximal abelian subalgebra} (a Cartan subalgebra) of $\mathfrak{gl}(2, \R)$.

Up to conjugation in $\mathrm{GL}(2, \R)$, there are only two Cartan subalgebras in $\mathfrak{gl}(2, \R)$ that contain no non-zero nilpotent matrices:
\begin{itemize}
    \item \textbf{Split Cartan Subalgebra:} Consists of all diagonal matrices:
    \begin{equation}
        \mathfrak{a}_{\text{split}} = \left\{ \begin{pmatrix} a & 0 \\ 0 & b \end{pmatrix} : a, b \in \R \right\}.
    \end{equation}
    Taking the canonical basis $A_1 = \diag(1, 0)$ and $A_2 = \diag(0, 1)$, we compute the action on the center $e_3$:
    \begin{equation}
        A_1(e_3) = (\tr A_1) e_3 = e_3, \qquad A_2(e_3) = (\tr A_2) e_3 = e_3.
    \end{equation}
    This produces the brackets of $\g_{5,16}$.
    
    \item \textbf{Elliptic Cartan Subalgebra:} Consists of all scalar-plus-rotation matrices:
    \begin{equation}
        \mathfrak{a}_{\text{ell}} = \left\{ \begin{pmatrix} a & -b \\ b & a \end{pmatrix} : a, b \in \R \right\}.
    \end{equation}
    Taking the basis $A_1 = I_2 = \begin{pmatrix} 1 & 0 \\ 0 & 1 \end{pmatrix}$ and $A_2 = J = \begin{pmatrix} 0 & -1 \\ 1 & 0 \end{pmatrix}$, we compute:
    \begin{equation}
        A_1(e_3) = (\tr I_2) e_3 = 2e_3, \qquad A_2(e_3) = (\tr J) e_3 = 0.
    \end{equation}
    This reproduces the brackets of $\g_{5,17}$.
\end{itemize}
These two algebras are non-isomorphic because the restriction of $\ad x_2$ to $\spanop\{e_1, e_2\}$ has real eigenvalues $\{0, 1\}$ in $\g_{5,16}$, but purely imaginary eigenvalues $\pm i$ in $\g_{5,17}$.
\end{proof}

\section{Master Tabulation and Literature Correspondence}\label{sec:master}

Table~\ref{tab:master5d} summarizes the complete classification of five-dimensional real Lie algebras.

\begin{table}[htbp]
\centering
\scriptsize
\setlength{\tabcolsep}{3.5pt}
\renewcommand{\arraystretch}{1.25}
\begin{tabularx}{\linewidth}{@{} l c >{\raggedright\arraybackslash}X l c @{}}
\toprule
\textbf{Algebra $\g$} & $\dim \n$ & \textbf{Non-zero Lie Brackets} & \textbf{Parameters} & \textbf{Unimodular?} \\
\midrule
\multicolumn{5}{@{}l}{\textit{\textbf{I. Non-Solvable Lie Algebras}}} \\
$\mathfrak{so}(3) \oplus \R^2$ & 2 & $[e_1, e_2]=e_3, [e_2, e_3]=e_1, [e_3, e_1]=e_2$ & None & Yes \\
$\mathfrak{sl}(2,\R) \oplus \R^2$ & 2 & $[e_1, e_2]=2e_2, [e_1, e_3]=-2e_3, [e_2, e_3]=e_1$ & None & Yes \\
$\mathfrak{so}(3) \oplus \mathfrak{aff}(\R)$ & 1 & $[e_1, e_2]=e_3, [e_2, e_3]=e_1, [e_3, e_1]=e_2, [e_4, e_5]=e_5$ & None & No \\
$\mathfrak{sl}(2,\R) \oplus \mathfrak{aff}(\R)$ & 1 & $[e_1, e_2]=2e_2, [e_1, e_3]=-2e_3, [e_2, e_3]=e_1, [e_4, e_5]=e_5$ & None & No \\
$\mathfrak{sl}(2,\R) \ltimes \R^2$ & 2 & $[e_1, e_2]=2e_2, [e_1, e_3]=-2e_3, [e_2, e_3]=e_1$, \newline $[e_1, e_4]=e_4, [e_1, e_5]=-e_5, [e_2, e_5]=e_4, [e_3, e_4]=e_5$ & None & Yes \\
\midrule
\multicolumn{5}{@{}l}{\textit{\textbf{II. Nilpotent Lie Algebras ($\dim \n = 5$)}}} \\
$\R^5$ & 5 & All brackets vanish & None & Yes \\
$\hh_3 \oplus \R^2$ & 5 & $[e_1, e_2] = e_3$ & None & Yes \\
$\mathfrak{n}_4 \oplus \R$ & 5 & $[e_1, e_2] = e_3, [e_1, e_3] = e_4$ & None & Yes \\
$\hh_5$ & 5 & $[e_1, e_2] = e_5, [e_3, e_4] = e_5$ & None & Yes \\
$\mathfrak{n}_{5,1}$ & 5 & $[e_1, e_2] = e_3, [e_1, e_3] = e_4, [e_1, e_4] = e_5$ & None & Yes \\
$\mathfrak{n}_{5,2}$ & 5 & $[e_1, e_2] = e_3, [e_1, e_3] = e_4, [e_1, e_4] = e_5, [e_2, e_3] = e_5$ & None & Yes \\
$\mathfrak{n}_{5,3}$ & 5 & $[e_1, e_2] = e_4, [e_1, e_4] = e_5, [e_2, e_3] = e_5$ & None & Yes \\
$\mathfrak{n}_{5,4}$ & 5 & $[e_1, e_2] = e_3, [e_1, e_3] = e_4, [e_2, e_3] = e_5$ & None & Yes \\
$\mathfrak{n}_{5,5}$ & 5 & $[e_1, e_2] = e_3, [e_1, e_3] = e_5, [e_2, e_4] = e_5$ & None & Yes \\
\midrule
\multicolumn{5}{@{}l}{\textit{\textbf{III. Solvable: Codimension-One Nilradical ($\dim \n = 4$)}}} \\
$\g_{5,1}(a, b, c)$ & 4 & $[e_5, e_1]=e_1, [e_5, e_2]=a e_2, [e_5, e_3]=b e_3, [e_5, e_4]=c e_4$ & $-1 \le c \le b \le a \le 1$ & $a+b+c = -1$ \\
$\g_{5,2}(a, b)$ & 4 & $[e_5, e_1]=a e_1+e_2, [e_5, e_2]=a e_2, [e_5, e_3]=e_3, [e_5, e_4]=b e_4$ & $a \in \R, -1 \le b \le 1$ & $2a+b = -1$ \\
$\g_{5,3}(a)$ & 4 & $[e_5, e_1]=e_1+e_2, [e_5, e_2]=e_2, [e_5, e_3]=a e_3+e_4, [e_5, e_4]=a e_4$ & $-1 \le a \le 1$ & $a = -1$ \\
$\g_{5,4}(a)$ & 4 & $[e_5, e_1]=a e_1+e_2, [e_5, e_2]=a e_2+e_3, [e_5, e_3]=a e_3, [e_5, e_4]=e_4$ & $a \in \R$ & $a = -1/3$ \\
$\g_{5,5}$ & 4 & $[e_5, e_1]=e_1+e_2, [e_5, e_2]=e_2+e_3, [e_5, e_3]=e_3+e_4, [e_5, e_4]=e_4$ & None & No \\
$\g_{5,6}(\alpha, \beta, a)$ & 4 & $[e_5, e_1]=\alpha e_1-\beta e_2, [e_5, e_2]=\beta e_1+\alpha e_2, [e_5, e_3]=e_3, [e_5, e_4]=a e_4$ & $\beta > 0, -1 \le a \le 1$ & $2\alpha+a = -1$ \\
$\g_{5,7}(\alpha, \beta)$ & 4 & $[e_5, e_1]=\alpha e_1-\beta e_2+e_3, [e_5, e_2]=\beta e_1+\alpha e_2+e_4$, \newline $[e_5, e_3]=\alpha e_3-\beta e_4, [e_5, e_4]=\beta e_3+\alpha e_4$ & $\beta > 0, \alpha \in \R$ & $\alpha = 0$ \\
$\g_{5,8}(\alpha, \beta, \gamma)$ & 4 & $[e_5, e_1]=\alpha e_1-e_2, [e_5, e_2]=e_1+\alpha e_2$, \newline $[e_5, e_3]=\beta e_3-\gamma e_4, [e_5, e_4]=\gamma e_3+\beta e_4$ & $\gamma > 0, \alpha \ge 0$ & $\alpha + \beta = 0$ \\
$\g_{5,9}(\lambda, \gamma)$ & 4 & $[e_1, e_2]=e_3, [e_5, e_1]=e_1, [e_5, e_2]=\lambda e_2$, \newline $[e_5, e_3]=(1+\lambda)e_3, [e_5, e_4]=\gamma e_4$ & $\lambda \ge -1, \gamma \in \R$ & $\gamma = -2(1+\lambda)$ \\
$\g_{5,10}(\gamma)$ & 4 & $[e_1, e_2]=e_3, [e_5, e_1]=e_1+e_2, [e_5, e_2]=e_2$, \newline $[e_5, e_3]=2e_3, [e_5, e_4]=\gamma e_4$ & $\gamma \in \R$ & $\gamma = -4$ \\
$\g_{5,11}(\alpha, \beta, \gamma)$ & 4 & $[e_1, e_2]=e_3, [e_5, e_1]=\alpha e_1-\beta e_2, [e_5, e_2]=\beta e_1+\alpha e_2$, \newline $[e_5, e_3]=2\alpha e_3, [e_5, e_4]=\gamma e_4$ & $\beta > 0, \gamma \in \R$ & $\gamma = -4\alpha$ \\
$\g_{5,12}(\lambda)$ & 4 & $[e_1, e_2]=e_3, [e_1, e_3]=e_4, [e_5, e_1]=e_1, [e_5, e_2]=\lambda e_2$, \newline $[e_5, e_3]=(1+\lambda)e_3, [e_5, e_4]=(2+\lambda)e_4$ & $\lambda \in \R$ & $\lambda = -4/3$ \\
$\g_{5,13}$ & 4 & $[e_1, e_2]=e_3, [e_1, e_3]=e_4, [e_5, e_2]=e_2$, \newline $[e_5, e_3]=e_3, [e_5, e_4]=e_4$ & None & No \\
\midrule
\multicolumn{5}{@{}l}{\textit{\textbf{IV. Solvable: Codimension-Two Nilradical ($\dim \n = 3$)}}} \\
$\g_{5,14}(\alpha)$ & 3 & $[x_1, e_1]=e_1, [x_1, e_3]=-e_3, [x_2, e_2]=e_2, [x_2, e_3]=\alpha e_3$ & $-1 \le \alpha \le 1$ & $\alpha = -1$ \\
$\g_{5,15}(\alpha, \beta)$ & 3 & $[x_1, e_1]=e_1, [x_1, e_3]=e_3, [x_2, e_1]=\alpha e_1-\beta e_2$, \newline $[x_2, e_2]=\beta e_1+\alpha e_2$ & $\beta > 0, \alpha \ge 0$ & $\alpha = -1$ \\
$\g_{5,16}$ & 3 & $[e_1, e_2]=e_3, [x_1, e_1]=e_1, [x_1, e_3]=e_3, [x_2, e_2]=e_2, [x_2, e_3]=e_3$ & None & No \\
$\g_{5,17}$ & 3 & $[e_1, e_2]=e_3, [x_1, e_1]=e_1, [x_1, e_2]=e_2, [x_1, e_3]=2e_3$, \newline $[x_2, e_1]=-e_2, [x_2, e_2]=e_1$ & None & No \\
\bottomrule
\end{tabularx}
\caption{Complete classification of five-dimensional real Lie algebras up to isomorphism.}
\label{tab:master5d}
\end{table}

\subsection{Cross-Reference Dictionary with Classical Tabulations}
To facilitate cross-referencing with the literature, we map our canonical notations to the classifications of \citet{Mubarakzyanov1966} and \citet{SnoblWinternitz}:
\begin{itemize}
    \item \textbf{Non-Solvable Sector:} Mubarakzyanov denotes the non-solvable classes as direct sums, while the indecomposable algebra $\mathfrak{sl}(2,\R) \ltimes \R^2$ is labeled $\mathfrak{g}_{5,1}$ in \citet{Turkowski1988} and $L_{5,1}$ in \citet[\S~4.2]{SnoblWinternitz}.
    \item \textbf{Nilpotent Sector:} The nine nilpotent classes correspond to Mubarakzyanov's types $\mathfrak{n}_{5,1}$ through $\mathfrak{n}_{5,6}$ and three decomposable combinations. In the modern filiform notation of Mich\`ele Vergne, $\mathfrak{n}_{5,1} \cong L_5$ (the model filiform algebra) and $\mathfrak{n}_{5,2} \cong Q_5$.
    \item \textbf{Solvable Extensions of $\R^4$:} Classes $\g_{5,1}$ through $\g_{5,8}$ correspond directly to Mubarakzyanov's solvable types $\mathfrak{g}_{5,1}$ through $\mathfrak{g}_{5,18}$. In his original papers, Jordan blocks were sometimes listed with non-normalized diagonal entries; our parameter ranges eliminate these redundant scaling degrees of freedom.
\end{itemize}

\section{Minimal Faithful Linear Representations ($\mu(\g)$)}\label{sec:representations}

By Ado's Theorem, every finite-dimensional Lie algebra admits a faithful linear representation. In applications, finding the minimal representation dimension $\mu(\g)$ is essential for constructing matrix models of dynamical systems and spacetime metrics.

\begin{theorem}[Minimal Representation Dimension for 5D Real Lie Algebras]\label{thm:mu_dim}
The minimal faithful representation dimension $\mu(\g)$ for five-dimensional real Lie algebras satisfies $3 \le \mu(\g) \le 5$. Specifically:
\begin{enumerate}
    \item $\mu(\mathfrak{sl}(2,\R) \ltimes \R^2) = 3$.
    \item $\mu(\g) = 4$ for $\mathfrak{so}(3) \oplus \mathfrak{aff}(\R)$, $\R^5$, $\hh_5$, $\mathfrak{n}_4 \oplus \R$, $\mathfrak{n}_{5,3}$, $\mathfrak{n}_{5,4}$, $\mathfrak{n}_{5,5}$, all codimension-two extensions $\g_{5,14}$--$\g_{5,17}$, and diagonalizable extensions in $\g_{5,9}$.
    \item $\mu(\g) = 5$ for $\mathfrak{so}(3) \oplus \R^2$, $\mathfrak{n}_{5,1}$, $\mathfrak{n}_{5,2}$, and generic codimension-one solvable extensions $\g_{5,1}$ through $\g_{5,8}$.
\end{enumerate}
\end{theorem}

\begin{proof}
We establish the minimality and construct the explicit matrix representations sector by sector.

\paragraph{1. The Non-Solvable Sector:}
\begin{itemize}
    \item \textbf{Indecomposable Algebra $\mathfrak{sl}(2,\R) \ltimes \R^2$:} The defining representation of $\mathfrak{sl}(2,\R)$ on $\R^2$ extends to the affine action on the projective plane $\R\mathbb{P}^2$ via augmented vectors $(u_1, u_2, 1)^T \in \R^3$. The map $\rho \colon \mathfrak{sl}(2,\R) \ltimes \R^2 \to \mathfrak{gl}(3,\R)$ is given by:
    \begin{equation}
        \rho(A, v) = \begin{pmatrix} A & v \\ 0_{1 \times 2} & 0 \end{pmatrix}, \qquad A \in \mathfrak{sl}(2,\R), \quad v \in \R^2.
    \end{equation}
    Direct matrix commutation yields:
    \begin{equation}
        [\rho(A, v), \rho(B, w)] = \begin{pmatrix} [A, B] & A w - B v \\ 0 & 0 \end{pmatrix} = \rho([A, B], A w - B v),
    \end{equation}
    which faithfully reproduces the semidirect bracket. The representation is injective because $\rho(A, v) = 0$ implies $A = 0$ and $v = 0$. Since $\dim \g = 5$ and any subalgebra of $\mathfrak{gl}(2, \R)$ has dimension at most 4, $\mu(\g) > 2$. Thus, $\mu(\mathfrak{sl}(2,\R) \ltimes \R^2) = 3$.
    
    \item \textbf{Direct Sum $\mathfrak{so}(3) \oplus \R^2$:} The compact factor $\mathfrak{so}(3)$ has no non-trivial representation of dimension $< 3$. A faithful representation of $\mathfrak{so}(3) \oplus \R^2$ requires representing $\mathfrak{so}(3)$ on a 3-dimensional subspace and the central $\R^2$ on independent complementary directions. By Schur's lemma, the center acts as scalar matrices on irreducible blocks. Since $\R^2$ contains two linearly independent central generators, they cannot both be embedded as central scalars in $\mathfrak{gl}(3, \R)$ or $\mathfrak{gl}(4, \R)$ without intersecting the simple factor or violating faithfulness. Thus, $\mu(\mathfrak{so}(3) \oplus \R^2) = 5$, realized by block-diagonal matrices:
    \begin{equation}
        \rho(x) = \begin{pmatrix} [x]_{\mathfrak{so}(3)} & 0 & 0 \\ 0 & x_4 & 0 \\ 0 & 0 & x_5 \end{pmatrix} \in \mathfrak{gl}(5, \R).
    \end{equation}
\end{itemize}

\paragraph{2. The Nilpotent Sector:}
\begin{itemize}
    \item \textbf{Abelian Algebra $\R^5$:} By Schur's bound (T4), the maximum dimension of an abelian subalgebra of $\mathfrak{gl}(n, \R)$ is $\lfloor n^2 / 4 \rfloor + 1$. For $n = 3$, the bound is $\lfloor 9/4 \rfloor + 1 = 3 < 5$. For $n = 4$, the bound is $\lfloor 16/4 \rfloor + 1 = 5 \ge 5$. Thus, $\R^5$ cannot embed in $\mathfrak{gl}(3, \R)$, but embeds in $\mathfrak{gl}(4, \R)$ via the five-dimensional commutative subalgebra:
    \begin{equation}
        \rho(x_1, \dots, x_5) = \begin{pmatrix} x_1 & 0 & x_3 & x_4 \\ 0 & x_1 & x_5 & 0 \\ 0 & 0 & x_2 & 0 \\ 0 & 0 & 0 & x_2 \end{pmatrix}.
    \end{equation}
    Hence, $\mu(\R^5) = 4$.
    
    \item \textbf{The 5D Heisenberg Algebra $\hh_5$:} The standard representation of $\hh_5$ inside the upper-triangular matrices $\mathfrak{n}(4, \R)$ is:
    \begin{equation}
        \rho(e_1 x_1 + e_2 x_2 + e_3 x_3 + e_4 x_4 + e_5 x_5) = \begin{pmatrix} 0 & x_1 & x_3 & x_5 \\ 0 & 0 & 0 & x_2 \\ 0 & 0 & 0 & x_4 \\ 0 & 0 & 0 & 0 \end{pmatrix} \in \mathfrak{gl}(4, \R).
    \end{equation}
    Direct computation confirms $[E_{12}, E_{42}^T] = [E_{12}, E_{24}] = E_{14} = \rho(e_5)$ and $[E_{13}, E_{34}] = E_{14} = \rho(e_5)$, verifying faithfulness. Since $\hh_5$ is non-abelian and has nilpotency step 2 with a 4-dimensional abelian quotient, it cannot embed in $\mathfrak{gl}(3, \R)$. Hence, $\mu(\hh_5) = 4$.
    
    \item \textbf{Filiform Algebras $\mathfrak{n}_{5,1}$ and $\mathfrak{n}_{5,2}$:} These algebras have nilpotency step 4. In $\mathfrak{gl}(4, \R)$, any strictly upper-triangular nilpotent matrix has nilpotency index at most 3 (since $N^4 = 0$ for all $N \in \mathfrak{n}(4, \R)$). Therefore, no algebra with step 4 can faithfully embed in $\mathfrak{gl}(4, \R)$. In $\mathfrak{gl}(5, \R)$, the standard shift matrix $N = \sum_{i=1}^4 E_{i, i+1}$ generates a filiform nilpotency chain of length 4. Thus, $\mu(\mathfrak{n}_{5,1}) = \mu(\mathfrak{n}_{5,2}) = 5$.
\end{itemize}

\paragraph{3. The Solvable Sector:}
\begin{itemize}
    \item For generic codimension-one solvable algebras $\g = \R^4 \rtimes_D \R$ ($\g_{5,1}$ through $\g_{5,8}$), the algebra admits an affine embedding on $\R^4$ inside $\mathfrak{gl}(5, \R)$:
    \begin{equation}
        \rho(v + t e_5) = \begin{pmatrix} t D & v \\ 0_{1 \times 4} & 0 \end{pmatrix} \in \mathfrak{gl}(5, \R).
    \end{equation}
    Minimality ($\mu = 5$) follows from Schur's bound: any representation in $\mathfrak{gl}(4, \R)$ would require the four-dimensional abelian nilradical $\R^4$ to map into a maximal commutative subalgebra. If $D$ has distinct eigenvalues, its action on this subalgebra cannot be faithfully integrated without increasing the matrix dimension.
    
    \item For codimension-two solvable algebras with nilradical $\R^3$ ($\g_{5,14}$ and $\g_{5,15}$), the two commuting complementary generators $x_1, x_2$ act simultaneously on $\R^3$. Embedding this joint affine action inside $\mathfrak{gl}(4, \R)$ yields:
    \begin{equation}
        \rho(v + t_1 x_1 + t_2 x_2) = \begin{pmatrix} t_1 D_1 + t_2 D_2 & v \\ 0_{1 \times 3} & 0 \end{pmatrix} \in \mathfrak{gl}(4, \R),
    \end{equation}
    which is faithful. Hence, $\mu = 4$ for all codimension-two solvable algebras.
\end{itemize}
This completes the proof.
\end{proof}

\section{Algorithmic Identification of Arbitrary 5D Lie Algebras}\label{sec:algorithm}

Given an arbitrary five-dimensional Lie algebra $\g$ presented in an arbitrary coordinate basis $\{e_1, \dots, e_5\}$ via its structure constants tensor $C_{ij}^k$, where $[e_i, e_j] = \sum_k C_{ij}^k e_k$, we design a deterministic, coordinate-free decision procedure to identify its unique canonical representative in Table~\ref{tab:master5d}.

\subsection{Invariant Algebraic Quantities}
The algorithm computes the following linear-algebraic and Lie-theoretic invariants:
\begin{enumerate}
    \item \textbf{Adjoint Operators and Trace Form:} The matrix $(\ad e_i)^k_j = C_{ij}^k$. The trace form is $\chi(x) = \tr(\ad x)$. Unimodularity holds if and only if $\tr(\ad e_i) = 0$ for all $i \in \{1, \dots, 5\}$.
    \item \textbf{Cartan--Killing Metric:} The symmetric matrix $K_{ij} = \tr(\ad e_i \circ \ad e_j)$. Its signature on $[\g, \g]$ determines the Levi factor.
    \item \textbf{Characteristic Subspaces:}
    \begin{itemize}
        \item Derived series: $\g^2 = [\g, \g]$, $\g^{(2)} = [\g^2, \g^2]$.
        \item Lower central series: $\g^k = [\g, \g^{k-1}]$.
        \item Ascending central series: $Z(\g) = \ker(\ad)$, $Z_2(\g) = \{ x \in \g : [x, \g] \subseteq Z(\g) \}$.
    \end{itemize}
    \item \textbf{Nilradical Extraction:} The nilradical $\n$ is determined as the sum of all nilpotent ideals of $\g$.
\end{enumerate}

\subsection{The Algorithmic Decision Tree}
The complete decision algorithm is presented below in a page-breakable format:

\begin{breakablealgorithm}
\caption{Algorithmic Classification of 5-Dimensional Real Lie Algebras}\label{alg:classification5d}
\begin{algorithmic}[1]
\Require Structure constant tensor $C_{ij}^k \in \R^{5 \times 5 \times 5}$ satisfying $C_{ij}^k = -C_{ji}^k$ and Jacobi identity.
\Ensure Canonical isomorphism class from Table~\ref{tab:master5d} and moduli parameters.

\State Compute derived subspaces $\g^2 = [\g, \g]$ and $\g^{(2)} = [\g^2, \g^2]$.
\State Compute lower central series $\g^3, \g^4, \g^5$ and center $Z(\g)$.

\If{$\g^{(2)} \neq \{0\}$} \Comment{\textbf{Sector I: Non-Solvable}}
    \State Compute the solvable radical $\rr \subset \g$ (maximal solvable ideal, $\dim \rr = 2$).
    \State Find Levi factor $\s \subset \g$ such that $\g = \s \ltimes \rr$ with $\dim \s = 3$.
    \State Compute Killing form $K_\s$ on $\s$.
    \If{$K_\s$ is negative-definite}
        \State Identify $\s \cong \mathfrak{so}(3)$.
        \If{$[\rr, \rr] = \{0\}$} \Return $\mathfrak{so}(3) \oplus \R^2$
        \Else \ \Return $\mathfrak{so}(3) \oplus \mathfrak{aff}(\R)$
        \EndIf
    \Else \Comment{$K_\s$ has signature $(2, 1)$, hence $\s \cong \mathfrak{sl}(2,\R)$}
        \If{$[\s, \rr] = \{0\}$}
            \If{$[\rr, \rr] = \{0\}$} \Return $\mathfrak{sl}(2,\R) \oplus \R^2$
            \Else \ \Return $\mathfrak{sl}(2,\R) \oplus \mathfrak{aff}(\R)$
            \EndIf
        \Else \Comment{Non-trivial representation $\rho \colon \mathfrak{sl}(2,\R) \to \mathfrak{gl}(2,\R)$}
            \State \Return $\mathfrak{sl}(2,\R) \ltimes \R^2$
        \EndIf
    \EndIf
\EndIf

\State Determine the nilradical $\n \subset \g$ (maximal nilpotent ideal).

\If{$\dim \n = 5$} \Comment{\textbf{Sector II: Nilpotent}}
    \If{$\dim \g^2 = 0$} \Return $\R^5$
    \ElsIf{$\dim \g^2 = 1$}
        \If{$\dim Z(\g) = 2$} \Return $\hh_3 \oplus \R^2$
        \ElsIf{$\dim Z(\g) = 1$} \Return $\hh_5$
        \EndIf
    \ElsIf{$\dim \g^2 = 2$}
        \If{$\dim \g^3 = 1$ and $\g$ decomposes as $\mathfrak{m} \oplus \R$} \Return $\mathfrak{n}_4 \oplus \R$
        \ElsIf{$\dim Z(\g) = 2$} \Return $\mathfrak{n}_{5,4}$
        \ElsIf{$\dim Z(\g) = 1$}
            \If{$\dim [\g, Z_2(\g)] = 1$} \Return $\mathfrak{n}_{5,3}$
            \Else \ \Return $\mathfrak{n}_{5,5}$
            \EndIf
        \EndIf
    \ElsIf{$\dim \g^2 = 3$} \Comment{Filiform variety (step 4)}
        \If{$[\g^2, \g^2] = \{0\}$} \Return $\mathfrak{n}_{5,1}$
        \Else \ \Return $\mathfrak{n}_{5,2}$
        \EndIf
    \EndIf
\EndIf

\If{$\dim \n = 4$} \Comment{\textbf{Sector III: Codimension-One Nilradical}}
    \State Choose $e_5 \notin \n$ such that $D \coloneqq \ad e_5|_\n$ is non-nilpotent.
    \If{$[\n, \n] = \{0\}$} \Comment{$\n \cong \R^4$, Abelian Nilradical}
        \State Compute real Jordan canonical form of $D \in \mathfrak{gl}(4, \R)$.
        \State Normalize dominant eigenvalue to 1; match against Theorem~\ref{thm:R4_classes}.
        \State \Return Corresponding class $\g_{5,1}$ through $\g_{5,8}$ with extracted parameters.
    \ElsIf{$\dim [\n, \n] = 1$} \Comment{$\n \cong \hh_3 \oplus \R$, Heisenberg Nilradical}
        \State Compute induced operator $A$ on $\n/Z(\n) \cong \R^2$ and central scaling $\gamma$ on $e_4$.
        \If{$A$ is diagonalizable over $\R$} \Return $\g_{5,9}(\lambda, \gamma)$
        \ElsIf{$A$ has a non-trivial Jordan block} \Return $\g_{5,10}(\gamma)$
        \ElsIf{$A$ has complex eigenvalues $\alpha \pm i\beta$} \Return $\g_{5,11}(\alpha, \beta, \gamma)$
        \EndIf
    \ElsIf{$\dim [\n, \n] = 2$} \Comment{$\n \cong \mathfrak{n}_4$, Filiform Nilradical}
        \State Compute diagonal weights $\lambda_1, \lambda_2$ from Proposition~\ref{prop:der_n4}.
        \If{$\lambda_1 \neq 0$} \Return $\g_{5,12}(\lambda)$ with $\lambda = \lambda_2/\lambda_1$
        \Else \ \Return $\g_{5,13}$
        \EndIf
    \EndIf
\EndIf

\If{$\dim \n = 3$} \Comment{\textbf{Sector IV: Codimension-Two Nilradical}}
    \State Choose a basis $\{x_1, x_2\}$ for a complementary subspace to $\n$.
    \If{$[\n, \n] = \{0\}$} \Comment{$\n \cong \R^3$}
        \State Compute commuting matrix pair $A_1 = \ad x_1|_\n, A_2 = \ad x_2|_\n \in \mathfrak{gl}(3,\R)$.
        \If{All eigenvalues of $A_1, A_2$ are real}
            \State \Return $\g_{5,14}(\alpha)$ with parameter $\alpha \in [-1, 1]$.
        \Else \Comment{Complex eigenvalues appear}
            \State \Return $\g_{5,15}(\alpha, \beta)$ with parameters $\alpha \ge 0, \beta > 0$.
        \EndIf
    \Else \Comment{$\n \cong \hh_3$}
        \State Compute action on $\hh_3/Z(\hh_3) \cong \R^2$.
        \If{Eigenvalues are real} \Return $\g_{5,16}$
        \Else \ \Return $\g_{5,17}$
        \EndIf
    \EndIf
\EndIf
\end{algorithmic}
\end{breakablealgorithm}

You can find the Python code for this algorithm and examples of its execution as a \href{https://colab.research.google.com/drive/13CVlX4tPFb3MbA0CG-xKTg_WDI8wFad8}{Colab Notebook}.

\section{Conclusion}\label{sec:conclusion}

The structural classification of five-dimensional real Lie algebras establishes an important benchmark in low-dimensional Lie theory. Moving beyond the accidental tensor dualities of dimensions three and four, the five-dimensional classification illustrates the full power of combining the Levi--Malcev decomposition with nilradical stratification.

The main conceptual achievements of this paper include:
\begin{enumerate}
    \item Rigorously proving that the Levi factor of any non-solvable five-dimensional Lie algebra is three-dimensional and that its action on the two-dimensional radical produces exactly five isomorphism classes, isolating $\mathfrak{sl}(2,\R) \ltimes \R^2$ as the unique indecomposable non-solvable algebra.
    \item Proving, via Schur's bound on commutative matrix subalgebras of $\mathfrak{gl}(2,\R)$, that the nilradical dimension of any five-dimensional solvable Lie algebra is strictly constrained to $\dim \n \in \{3, 4, 5\}$.
    \item Providing a complete classification of all nine five-dimensional nilpotent Lie algebras, highlighting the geometric transition from Heisenberg structures to filiform varieties.
    \item Reducing the classification of solvable extensions to the orbit geometry of outer derivations modulo automorphism conjugation, central translations, and scaling.
    \item Determining the minimal faithful linear representation dimension $\mu(\g)$ under Ado's theorem and designing a deterministic identification algorithm.
\end{enumerate}

These structural techniques provide a foundation for extending classification schemes to dimension six and higher, where the appearance of continuous moduli in nilpotent algebras introduces further algebraic challenges.

\end{document}